\documentclass[11pt,reqno]{amsart}
\usepackage{amssymb,mathrsfs,graphicx,subfigure, enumerate}
\usepackage{amsmath,amsfonts,amssymb,amscd,amsthm,bbm}
\usepackage{extpfeil}
\usepackage{graphicx,colortbl}
\usepackage{float}
\usepackage{epsfig}
\usepackage{caption}
\usepackage{bm}
\graphicspath{{Figure/}}
\usepackage{tikz}
\usetikzlibrary{arrows.meta,positioning,calc}
\usepackage[margin=1in]{geometry}

\makeatletter
\@namedef{subjclassname@2020}{\textup{2020} Mathematics Subject Classification}
\makeatother

\title[Incompressible Schr\"odinger flow]{Local well-posedness of the two-dimensional incompressible Schr\"odinger flow}

\author[Kim]{Jeongho Kim}
\address[Jeongho Kim]{\newline Department of Applied Mathematics, \newline 
	Kyung Hee University, 1732, Deogyeong-daero, Giheung-gu, \newline
	Yongin-si, Gyeonggi-do 17104, Republic of Korea}
\email{jeonghokim@khu.ac.kr}

\newtheorem{theorem}{Theorem}[section]
\newtheorem{lemma}{Lemma}[section]

\newtheorem{proposition}{Proposition}[section]

\newtheorem{remark}{Remark}[section]

\newtheorem{definition}{Definition}[section]

\newcommand{\bbr}{\mathbb R}
\newcommand{\bbc}{\mathbb C}
\newcommand{\bbs}{\mathbb S}

\newcommand{\bbn} {\mathbb N}

\newcommand{\<}{\langle}
\renewcommand{\>}{\rangle}
\newcommand{\tr}{\textup{tr}}
\renewcommand{\Im}{\textup{Im}}
\renewcommand{\Re}{\textup{Re}}

\renewcommand{\i}{\textup{i}}
\newcommand{\dg}{\dagger}

\newcommand{\e}{\varepsilon}
\newcommand{\pa}{\partial}
\newcommand{\dv}{\textup{div}}
\newcommand{\curl}{\textup{curl}}
\begin{document}
%%%%%%%%%%%%%%%%


\thanks{\textbf{Acknowledgment.} J. Kim was supported by Samsung Science and Technology Foundation under Project Number SSTF-BA2401-01 }

\begin{abstract}
	We consider the incompressible Schr\"odinger flow (ISF), an efficient  model for simulating incompressible inviscid fluids introduced by Chern et al. \cite{C17,CKPSW16}, on a bounded, simply connected two-dimensional domain. After transforming it into the director formulation, the equation becomes the Schr\"odinger map (or the Landau--Lifshitz equation) transported by a velocity field generated from the director itself via Biot--Savart law. We show that the director formulation has a unique local-in-time strong solution. Then, we lift the solution to the director formulation to the spinor, and conclude that it is indeed a unique solution to the ISF. This result is a first rigorous analysis on ISF, a widely used computational model.
\end{abstract}
\maketitle 
%\tableofcontents

\section{Introduction}

We consider the incompressible Schr\"odinger flow (ISF) on a two-dimensional bounded, simply connected domain $\Omega\subset\bbr^2$ with smooth boundary $\pa\Omega$. Specifically, let $\psi=(\psi_1,\psi_2)^\top:\Omega\to \bbc^2$ be a {\it spinor} field. Then, the governing equation for the ISF reads as follows:
\begin{equation}\label{eq:ISF}
	\i\hbar\pa_t \psi +\frac{\hbar^2}{2}\Delta\psi = p\psi, \quad t>0,\quad x\in\Omega,
\end{equation}
subject to the initial data $\psi_0=(\psi_{01},\psi_{02})$ satisfying $|\psi_0|^2=1$ and the boundary condition $\pa_\nu \psi=0$ on $\pa\Omega$, where $\nu$ denotes the outward unit normal vector of the domain $\Omega$.
Here, $\hbar$ is a model parameter, which is reminiscent of the Planck constant in the Schr\"odinger equation, and $p$ is a real Lagrange multiplier that enforces the following constraint on the spinor field:
\begin{equation}\label{eq:incompressibility}
	\Im(\psi^\dg\Delta\psi)=0,\quad \psi^\dg:=(\overline{\psi_1},\overline{\psi_2}).
\end{equation}
Indeed, if we define the ``velocity" $v:\Omega\to\bbr^2$ corresponds to $\psi$ as
\[v=(v_1,v_2)=(\hbar\Im(\psi^\dg\pa_1\psi),\hbar\Im(\psi^\dg\pa_2\psi)),\]
the constraint \eqref{eq:incompressibility} can be represented as the incompressibility condition of $v$:
\[\dv\,v = \hbar\Im(\psi^\dg\Delta\psi)=0.\]
Furthermore, the condition \eqref{eq:incompressibility} yields the additional constraint $|\psi|^2=\psi^\dg\psi=1$, since the initial data satisfies $|\psi_0|^2=\psi^\dg_0\psi_0 =1$, and for all $t>0$,
\begin{align*}
	\frac{d}{dt}(\psi^\dg\psi)&=\left(\frac{\i\hbar}{2}\Delta\psi-\frac{\i}{\hbar}p\psi\right)^\dg\psi + \psi^\dg\left(\frac{\i\hbar}{2}\Delta\psi-\frac{\i}{\hbar}p\psi\right)\\
	&=\frac{\i\hbar}{2}\left(\psi^\dg\Delta\psi-(\Delta\psi)^\dg\psi\right)=-\hbar\Im(\psi^\dg\Delta\psi)=0.
\end{align*}
Therefore, the spinor field $\psi\in \bbc^2$ with unit length can be identified as a vector in $\bbs^3\subset\bbc^2$. 

The ISF \eqref{eq:ISF} was first introduced by Chern et al. \cite{CKPSW16}, as an effective model for simulating incompressible Euler equations, which enables robust simulation of the fine structure of the vortices. See also \cite{C17,CKPS17,YXZFLZ21}
for subsequent developments on this direction. To identify the relation between the ISF \eqref{eq:ISF} and the incompressible Euler equations, we define $n=\psi^\dg\sigma\psi$, called the {\it director} field, borrowing the term from the theory of liquid crystals \cite{GP93,LL95}. Here, 
\[\sigma := (\sigma_1,\sigma_2,\sigma_3),\quad\sigma_1 = \begin{pmatrix}
	0 & 1 \\ 1 & 0
\end{pmatrix},\quad \sigma_2 = \begin{pmatrix}
	0 & -\i \\ \i & 0
\end{pmatrix},\quad \sigma_3 = \begin{pmatrix}
	1 & 0 \\ 0 & -1
\end{pmatrix},\]
are the Pauli matrices. In particular, if the spinor field $\psi\in\bbs^3$ has unit length, then it is straightforward to verify that the director field $n=(n_1,n_2,n_3)\in \bbr^3$ also has a unit length, identifying it with the $\bbs^2$-valued vector field. Moreover, it follows from the ISF \eqref{eq:ISF} that the director $n$ satisfies the following equation:

\begin{equation}\label{eq:ISF_director}
	\pa_t n +(v\cdot\nabla)n = \frac{\hbar}{2}n\times\Delta n,
\end{equation}
subject to the initial data $n_0=\psi_0^\dg\sigma\psi_0$ and the boundary condition $\pa_\nu n=0$ on $\pa\Omega$. We refer to Section \ref{sec:2} for the detailed derivation of \eqref{eq:ISF_director} from \eqref{eq:ISF}. In particular, the velocity field $v=\hbar\Im(\psi^\dg\nabla\psi)$ can be completely represented in terms of the director $n$ via the Biot--Savart law:
\begin{equation}\label{v:Kn}
	v=K[n]:=\nabla^\perp \Delta^{-1}_D\omega,\quad \omega:=\frac{\hbar}{2}n\cdot(\pa_1n\times\pa_2 n).
\end{equation}
Here, $\nabla^\perp:=(-\pa_2,\pa_1)$, and $\Delta_D^{-1}f$ denotes the solution of the Poisson equation $\Delta u = f$ on $\Omega$ with a homogeneous Dirichlet boundary condition $u\equiv0$ on $\pa\Omega$. In particular, the vorticity $\omega$ follows the two-dimensional Mermin--Ho relation of superfluid \cite{MH76}. 

Furthermore, the velocity field $v$ satisfies its own time evolution equation, which is given by
\begin{equation}\label{eq:ISF_velocity}
	\pa_t v +(v\cdot\nabla) v +\nabla p =-\frac{\hbar^2}{4}\nabla\cdot(\nabla n\odot \nabla n),\quad \dv\,v=0,
\end{equation}
subject to the initial data $v_0 = \hbar\Im(\psi^\dg_0\nabla\psi_0)$ and the slip boundary condition $v\cdot\nu=0$, where $(\nabla n\odot \nabla n)_{jk}:=\pa_j n \cdot\pa_k n$ so that the right-hand side of \eqref{eq:ISF_velocity} is the Ericksen stress familiar from the hydrodynamics of liquid crystals \cite{LL95}. Again, we refer to Section \ref{sec:2} for the derivation of \eqref{eq:ISF_velocity}. As a consequence, when $\hbar\to0$, we expect that the velocity field $v$ formally converges to the incompressible inviscid flow, and this is the motivation of using the ISF as an effective simulation model for the incompressible Euler equations. This shares the same spirit of the hydrodynamic formulation of Schr\"odinger-type equations, which goes back to Madelung \cite{Ma27}.

Interestingly, the director formulation \eqref{eq:ISF_director} can be understood as a generalization of so--called the Schr\"odinger map (or Schr\"odinger flow). Indeed, when the velocity $v$ vanishes, the equation \eqref{eq:ISF_director} becomes the Schr\"odinger map
\begin{equation}\label{eq:LL}
	\pa_t n =n\times\Delta n,
\end{equation}
which is also called the Landau--Lifshitz equation for the Heisenberg model for a ferromagnetic spin \cite{LL35,GD08}. Over decades, \eqref{eq:LL} and related models have been extensively studied in literature. To name a few, the global existence of weak solution, local well-posedness of the smooth solution of \eqref{eq:LL} in the whole space was studied in \cite{SSB86}. A generalization of the Schr\"odinger map, whose image is a general compact Riemann surface was considered in \cite{CSU00}, and a Schr\"odinger map from a compact Riemannian manifold to a symplectic or K\"ahler manifold was investigated in \cite{DW98, DW01, M07}. When the domain is bounded, Chen and Wang \cite{CW23} proved the existence and uniqueness of the local strong solution to the initial-boundary value problem for \eqref{eq:LL}, where the parabolic perturbation approximation considered in \cite{CF01,CJ18} was crucially used. When the flow is also convected by the prescribed incompressible velocity $v$, that is, when $v$ in \eqref{eq:ISF_director} is a given divergence-free vector field, the global existence of weak solution, as well as the local existence of strong solution was studied in \cite{CWW26}. We also refer to \cite{HZ24}, where the so-called Navier--Stokes--Schr\"odinger system, in which \eqref{eq:ISF_director} coupled with the viscous version of \eqref{eq:ISF_velocity}, was investigated. For more comprehensive and extensive literature regarding the Schr\"odinger map or Landau--Lifshitz equations, we refer to \cite{CW23,CWW26} and references therein. However, to the extent of the author's knowledge, the system \eqref{eq:ISF_director}, where the velocity is determined as a nonlinear function of the director, has not been treated in the previous work, and it is interesting to establish the well-posedness theory on it, since it is directly related to the theoretical justification of the ISF for simulating incompressible Euler equations.

Therefore, the first goal of this paper is to establish local well-posedness of \eqref{eq:ISF_director}. Basically, we follow the argument introduced in \cite{CWW26}, that is, introducing the regularized system \eqref{eq:ISF_regularize} and then taking the limit $\e\to0$. The difficulty of \eqref{eq:ISF_director} compared with the prescribed velocity case considered in \cite{CWW26} is that the velocity itself is a nonlinear function of the director $n$. To overcome this nonlinearity, we construct the sequence of approximated solutions by using the iterative method, and carefully track the dependency of $v$ on $n$, obtaining the uniform estimate during the iteration. Additionally, we attain the uniqueness of the solution, which is not considered in \cite{CWW26}. Before stating the first main result, we present the definition of solutions.

\begin{definition}\label{def:director}
	Let $T>0$ and $n_0\in H^3(\Omega;\bbs^2)$. We say that $n$ is a strong solution to \eqref{eq:ISF_director}--\eqref{v:Kn} on $[0,T]$ with the initial data $n_0$ if
	\[n\in L^\infty(0,T;H^3(\Omega;\bbs^2)),\quad \pa_t n\in L^\infty(0,T;H^1(\Omega)),\]
	the equation \eqref{eq:ISF_director} holds almost everywhere in $[0,T]\times\Omega$ with $v=K[n]$, $\pa_\nu n=0$ on $\pa\Omega$ in the sense of traces for $t\in(0,T)$, and $n(0)=n_0$.
\end{definition}

Then, the first main result is the following.

\begin{theorem}\label{thm:director}
	Let $n_0\in H^3(\Omega;\bbs^2)$ with $\pa_\nu n_0=0$ on $\pa\Omega$. Then, there exist positive constants $T_*$ and $M$, depending only on the initial data $n_0$ and $\hbar$, such that there exists a unique strong solution $n$ to \eqref{eq:ISF_director} on $[0,T_*]$ in the sense of Definition \ref{def:director} satisfying 
	\[\sup_{0\le t\le T_*}\left(\|n(t)\|_{H^3}^2+\|\pa_tn(t)\|_{H^1}^2\right)\le M.\]
	Moreover, $n\in C([0,T_*];H^{3-\delta})$ for every $\delta>0$.
\end{theorem}

The second main result concerns the well-posedness of the original ISF equation \eqref{eq:ISF}. Once we get the strong solution $n$ to \eqref{eq:ISF_director} from Theorem \ref{thm:director}, we use it to construct the spinor $\psi$ that satisfies $n=\psi^\dg\sigma \psi$, the incompressibility condition $\Im(\psi^\dg\Delta\psi)=0$, and $v=\hbar\Im(\psi^\dg\nabla\psi)$. Then, we directly show that this constructed spinor $\psi$ indeed satisfies the ISF equation \eqref{eq:ISF}. Again, we first fix the notion of solution.

\begin{definition}\label{def:spinor}
	Let $T>0$ and $\psi_0\in H^3(\Omega;\bbs^3)$. We say that $(\psi,p)$ is a strong solution to \eqref{eq:ISF} on $[0,T]$ with the initial data $\psi_0$ if
	\[\psi\in L^\infty(0,T;H^3(\Omega;\bbs^3)),\quad \pa_t\psi\in L^\infty(0,T;H^1(\Omega)),\quad p\in L^\infty(0,T;L^2(\Omega))\ \mbox{ is real-valued},\]
	the equation \eqref{eq:ISF} and the constraint \eqref{eq:incompressibility} hold almost everywhere in $[0,T]\times\Omega$, $\pa_\nu\psi=0$ on $\pa\Omega$ in the sense of traces for $t\in(0,T)$, and $\psi(0)=\psi_0$.
\end{definition}

Then, the precise statement of the well-posedness of ISF reads as follows.

\begin{theorem}\label{thm:spinor}
	Let $\psi_0\in H^3(\Omega;\bbs^3)$ satisfy $\pa_\nu\psi_0=0$ on $\pa\Omega$ and $\dv\,\Im(\psi^\dg_0\nabla\psi_0)=0$ in $\Omega$. Then, there exist positive constants $T_*$ and $M$, depending only on the initial data $\psi_0$ and $\hbar$, such that there exists a unique strong solution $(\psi,p)$ to \eqref{eq:ISF} on $[0,T_*]$ in the sense of Definition \ref{def:spinor} satisfying the normalization
	\[\int_{\Omega}p(t,x)\,dx=0\quad\mbox{for}\quad t\in[0,T_*]\]
	and
	\[\sup_{0\le t\le T_*}\left(\|\psi(t)\|_{H^3}^2+\|\pa_t\psi(t)\|_{H^1}^2\right)\le M.\]
	Moreover, $p$ has better regularity $p\in L^\infty(0,T_*;H^1)$.
	
\end{theorem}

\begin{remark}\label{rem:uniqueness}
	The normalization condition $\int_\Omega p(t,x)\,dx=0$ is essential for the uniqueness of the solution. Indeed, for any Lipschitz continuous function $\theta:[0,T_*]\to\bbr$ with $\theta(0)=0$, if $(\psi,p)$ is solution then,
	\[(\tilde{\psi},\tilde{p}) := (e^{\i\theta(t)}\psi,p-\hbar\theta'(t))\] 
	is also a solution to \eqref{eq:ISF}. Therefore, the solution is uniquely determined after normalizing the pressure as $\int_\Omega p\,dx=0$. Under this restriction, we have $\theta'(t)=0$ for $t\in[0,T_*]$, which implies $\theta(t)\equiv0$, and $(\psi,p)\equiv(\tilde{\psi},\tilde{p})$.
\end{remark}

Figure \ref{fig:schematic} below shows a brief relation between the spinor equation \eqref{eq:ISF} and the director equation \eqref{eq:ISF_director}. In particular, the solution $\psi$ \eqref{eq:ISF} is constructed from the solution $n$ to \eqref{eq:ISF_director} by lifting the solution via unitary matrix $U(n,n_0)$ and the phase correction $\theta$ that satisfies the Poisson equation.

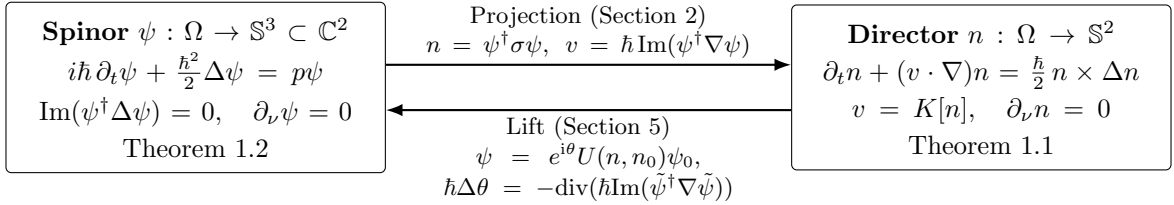
\begin{figure}[h!]
\begin{tikzpicture}[
	box/.style={draw, rounded corners=2pt, align=center, inner sep=6pt,
		minimum height=1.6cm, text width=4.6cm, font=\small},
	arr/.style={-{Latex[length=2.2mm]}, thick},
	lab/.style={font=\footnotesize, align=center, inner sep=2pt}
	]
	\node[box] (psi) at (0,0) {\textbf{Spinor} $\psi:\Omega\to\mathbb S^3\subset\mathbb C^2$\\[2pt]
		$i\hbar\,\partial_t\psi+\frac{\hbar^2}{2}\Delta\psi=p\psi$\\[2pt]
		$\mathrm{Im}(\psi^\dagger\Delta\psi)=0,\quad \partial_\nu\psi=0$\\[2pt] {Theorem 1.2}};
	\node[box] (n) at (10.4,0) {\textbf{Director} $n:\Omega\to\mathbb S^2$\\[2pt]
		$\partial_tn+(v\cdot\nabla)n=\frac{\hbar}{2}\,n\times\Delta n$\\[2pt]
		$v=K[n],\quad \partial_\nu n=0$\\[2pt] {Theorem 1.1}};
	
	\draw[arr] ([yshift=3mm]psi.east) -- ([yshift=3mm]n.west)
	node[lab, midway, above, text width=5.2cm] {Projection (Section 2)\\ $n=\psi^\dagger\sigma\psi$,\ \ $v=\hbar\,\mathrm{Im}(\psi^\dagger\nabla\psi)$};
	\draw[arr] ([yshift=-3mm]n.west) -- ([yshift=-3mm]psi.east)
	node[lab, midway, below, text width=5.2cm] {Lift (Section 5)\\ $\psi=e^{\i\theta}U(n,n_0)\psi_0$,\ \ $\hbar\Delta\theta=-\dv(\hbar\Im(\tilde{\psi}^\dg\nabla\tilde{\psi}))$};
\end{tikzpicture}
\caption{Schematic relation between the spinor equation (ISF) \eqref{eq:ISF} and the director equation \eqref{eq:ISF_director} and their solutions.}
\label{fig:schematic}
\end{figure}

The rest of the paper is organized as follows. In Section \ref{sec:2}, we formally derive the equations \eqref{eq:ISF_director} for the director $n$ and \eqref{eq:ISF_velocity} for the velocity $v$ from the ISF equation \eqref{eq:ISF}, justifying the equivalence between these systems. Section \ref{sec:3} introduces the $\e$-regularized system \eqref{eq:ISF_regularize} and prove the existence of the strong solution by using iteration method with uniform estimate. Then, in Section \ref{sec:4}, we obtain the uniform-in-$\e$ estimate for the solution to the regularized system to show that the existence time is indeed $\e$-independent. Then, we show that the limit as $\e\to0$ is indeed the unique strong solution to \eqref{eq:ISF_director}. Finally, in Section \ref{sec:5}, we construct the unique strong solution $(\psi,p)$ to the ISF \eqref{eq:ISF} by lifting the strong solution $n$ to \eqref{eq:ISF_director}.\\

\section{Derivation of \eqref{eq:ISF_director} and \eqref{eq:ISF_velocity}}\label{sec:2}
\setcounter{equation}{0}

In this section, we derive the equations \eqref{eq:ISF_director} and \eqref{eq:ISF_velocity} from the ISF equation \eqref{eq:ISF}. Before we present the detailed derivation, we first introduce several identities regarding the Pauli matrices, which will be used throughout the paper. First of all, for any $\psi\in \bbc^2$ with $|\psi|^2=1$, and $n=\psi^\dg \sigma\psi$, we have
\begin{equation}\label{projector}
	\psi\psi^\dg = \frac{1}{2}(I+n\cdot\sigma).
\end{equation}
Furthermore, it is well-known that the Pauli's matrices satisfy the following identities:
\[\sum_{k=1}^3(\sigma_k)_{ab}(\sigma_k)_{cd} = 2\delta_{ad}\delta_{bc}-\delta_{ab}\delta_{cd},\quad\mbox{and}\quad  \sigma_j\sigma_k = \delta_{jk}I+\i\e_{jkl}\sigma_l.\]
Therefore, for any $\phi,\psi\in \bbc^2$, we get
\begin{equation}\label{Pauli_identity_1}
	\sum_{k=1}^3(\psi^\dg\sigma_k\phi)\sigma_k = 2\phi\psi^\dg -(\psi^\dg\phi)I,\quad\mbox{hence}\quad \phi=\frac{1}{2}\left((\psi^\dg\phi)\psi+\left((\psi^\dg\sigma\phi)\cdot\sigma\right)\psi\right),
\end{equation}
and for any $a,b\in\bbc^3$,
\begin{equation}\label{Pauli_identity_2}
	(a\cdot\sigma)(b\cdot\sigma) = \sum_{j,k}a_j\sigma_jb_k\sigma_k =\sum_{j,k}a_jb_k(\delta_{jk}I+\i\e_{jkl}\sigma_l)=(a\cdot b)I+\i(a\times b)\cdot\sigma.
\end{equation}

\subsection{Derivation of \eqref{eq:ISF_director} from \eqref{eq:ISF}}
We first derive the director equation \eqref{eq:ISF_director} from the ISF \eqref{eq:ISF}. Using \eqref{eq:ISF} and the definition $n=\psi^\dg\sigma\psi$, we compute 
\begin{align}
	\begin{aligned}\label{patn}
	\pa_t n &= (\pa_t\psi)^\dg\sigma\psi+\psi^\dg\sigma(\pa_t\psi)\\
	&=\left(\frac{\i\hbar}{2}\Delta\psi-\frac{\i}{\hbar}p\psi\right)^\dg \sigma\psi +\psi^\dg\sigma\left(\frac{\i\hbar}{2}\Delta\psi -\frac{\i}{\hbar}p\psi\right)\\
	&=\frac{\i\hbar}{2}(\psi^\dg\sigma\Delta \psi-(\Delta \psi)^\dg\sigma\psi)=-\hbar\Im(\psi^\dg\sigma\Delta\psi)=-\hbar\pa_j\Im(\psi^\dg\sigma\pa_j\psi).
	\end{aligned}
\end{align}
For fixed $j\in\{1,2\}$, define $u_j:=\psi^\dg\sigma\pa_j\psi\in\bbc^3$. Then,
\begin{equation}\label{Rew}
	\Re(u_j) = \Re(\psi^\dg \sigma\pa_j\psi)=\frac{1}{2}\pa_j(\psi^\dg\sigma\psi)=\frac{1}{2}\pa_j n.
\end{equation}
On the other hand, applying the identity \eqref{Pauli_identity_1} with  $\phi=\pa_j\psi$, we get
\begin{equation}\label{eq:pa_jpsi}
	\pa_j\psi = \frac{1}{2}\left((\psi^\dg\pa_j\psi)\psi+\left(u_j\cdot\sigma\right)\psi\right)=\frac{1}{2}\left(\frac{\i v_j}{\hbar}\psi +(u_j\cdot\sigma)\psi\right),
\end{equation}
where the last equality holds since $|\psi|^2=1$ and
\begin{equation}\label{psidgpajpsi}
	\psi^\dg\pa_j\psi = \Re(\psi^\dg\pa_j\psi)+\i\Im(\psi^\dg\pa_j\psi)=\frac{1}{2}\pa_j|\psi|^2+\frac{\i v_j}{\hbar}=\frac{\i v_j}{\hbar},
\end{equation}
using the definition $v=\hbar\Im(\psi^\dg\nabla\psi)$. Therefore, multiplying $\psi^\dg\sigma_k$ on the left of \eqref{eq:pa_jpsi} and using \eqref{Pauli_identity_2} with $a=e_k$ and $b=u_j$, we get
\begin{align*}
	(u_j)_k &=\psi^\dg\sigma_k\pa_j\psi= \frac{\i v_j}{2\hbar}n_k + \frac{1}{2}\psi^\dg\sigma_k(u_j\cdot\sigma)\psi=\frac{\i v_j}{2\hbar} n_k+\frac{1}{2}\psi^\dg((u_j)_k I +\i(e_k\times u_j)\cdot\sigma)\psi\\
	&=\frac{\i v_j}{2\hbar}n_k+\frac{(u_j)_k}{2}+\frac{\i}{2}(u_j\times n)_k.
\end{align*}
where we used
\[\psi^\dg(e_k\times u_j)\cdot\sigma \psi = (e_k\times u_j)\cdot n = (u_j\times n)_k.\]
Therefore, $u_j$ satisfies the following identity
\[u_j = \frac{\i v_j}{\hbar} n +\i u_j\times n,\]
whose imaginary part can be computed by using \eqref{Rew} as
\begin{equation}\label{Imw}
	\Im(u_j) = \frac{v_j}{\hbar}n + \Re(u_j)\times n=\frac{v_j}{\hbar}n+\frac{1}{2}\pa_jn\times n.
\end{equation}
Returning to \eqref{patn}, we use \eqref{Imw} to observe that $n$ satisfies
\[\pa_t n = -\hbar\pa_j \Im(u_j)=-\pa_j(v_j n)-\frac{\hbar}{2}\pa_j(\pa_j n\times n)=-(v\cdot\nabla) n +\frac{\hbar}{2}n\times\Delta n,\]
since $\dv\,v=\hbar\Im(\psi^\dg\Delta \psi)=0$. Next, we represent the velocity $v$ in terms of $n$. We first observe that
\[\curl\,v = \pa_1 v_2-\pa_2 v_1 = 2\hbar\Im((\pa_1\psi)^\dg\pa_2\psi),\quad \dv\,v=0.\]
Since $\Omega$ is simply connected, $v=\nabla^\perp\phi$ for some scalar function $\phi$, which satisfies $\Delta\phi = \curl\,v=2\hbar\Im((\pa_1\psi)^\dg\pa_2\psi)$. Moreover, we observe that for the tangential vector $\nu^\perp=(-\nu_2,\nu_1)$ of $\pa\Omega$, 
\[\pa_{\nu^\perp}\phi = -\nu_2\pa_1\phi+\nu_1\pa_2\phi =-\nabla^\perp\phi \cdot\nu =-v\cdot\nu = -\hbar\Im(\psi^\dg\pa_\nu\psi)=0. \]
This implies that $\phi$ is constant along the boundary, and it can be normalized to be 0. Hence, $v$ can be represented as
\[v =\nabla^\perp\Delta_D^{-1}(2\hbar\Im((\pa_1\psi)^\dg\pa_2\psi)).\]
Therefore, it suffice to represent the quantity $2\hbar\Im((\pa_1\psi)^\dg\pa_2\psi)$ in terms of $n$. To this end, let $P:=\psi\psi^\dg$ be a projector defined as in \eqref{projector} and we compute $T:=\tr(P[\pa_1P,\pa_2P])$ in two different ways, where $[A,B]:=AB-BA$. 

First, it follows from \eqref{projector} that $\pa_jP = \frac{1}{2}\pa_j n\cdot \sigma$. Consequently, using \eqref{Pauli_identity_2}, we get
\[[\pa_1P,\pa_2P] = \frac{1}{4}[\pa_1 n\cdot\sigma,\pa_2n\cdot\sigma]=\frac{\i}{2}(\pa_1n\times\pa_2 n)\cdot\sigma,\]
and we use $P=\frac{1}{2}(I+n\cdot\sigma)$ once again to get
\begin{equation}\label{represent1}
	\tr(P[\pa_1P,\pa_2P]) = \frac{\i}{4}\tr\left((I+n\cdot\sigma)(\pa_1 n\times\pa_2 n)\cdot\sigma\right)=\frac{\i}{2}n\cdot(\pa_1n\times\pa_2n),
\end{equation}
since $\tr\,\sigma=0$. 

On the other hand, we compute $\tr(P[\pa_1P,\pa_2P])$ in another way as
\[\tr(P[\pa_1P,\pa_2P]) = \tr(\psi^\dg[\pa_1P,\pa_2P]\psi)=\psi^\dg[\pa_1P,\pa_2P]\psi.\]
Since $\pa_jP = (\pa_j\psi)\psi^\dg +\psi(\pa_j\psi)^\dg$, we have
\begin{align*}
\psi^\dg(\pa_1P)(\pa_2P)\psi &= \psi^\dg((\pa_1\psi) \psi^\dg +\psi(\pa_1\psi)^\dg)((\pa_2\psi) \psi^\dg +\psi(\pa_2\psi)^\dg)\psi\\
&=(\psi^\dg(\pa_1\psi)\psi^\dg +(\pa_1\psi)^\dg)(\pa_2\psi+\psi(\pa_2\psi)^\dg\psi)\\
&=\psi^\dg (\pa_1\psi)\psi^\dg(\pa_2\psi)+(\pa_1\psi)^\dg(\pa_2\psi)+\psi^\dg(\pa_1\psi)(\pa_2\psi)^\dg\psi+(\pa_1\psi)^\dg\psi(\pa_2\psi)^\dg\psi.
\end{align*}
However, it follows from \eqref{psidgpajpsi} that
\begin{align*}
	&\psi^\dg (\pa_1\psi)\psi^\dg(\pa_2\psi)= \left(\frac{\i v_1}{\hbar}\right)\left(\frac{\i v_2}{\hbar}\right)=-\frac{v_1v_2}{\hbar^2},\\
	&\psi^\dg(\pa_1\psi)(\pa_2\psi)^\dg\psi = \left(\frac{\i v_1}{\hbar}\right)\left(-\frac{\i v_2}{\hbar}\right)=\frac{v_1v_2}{\hbar},\\
	&(\pa_1\psi)^\dg\psi(\pa_2\psi)^\dg\psi=\left(-\frac{\i v_1}{\hbar}\right)\left(-\frac{\i v_2}{\hbar}\right)=-\frac{v_1v_2}{\hbar^2},
\end{align*}
which yields
\begin{equation}\label{represent2}
	\psi^\dg[\pa_1P,\pa_2P]\psi = (\pa_1\psi)^\dg(\pa_2\psi)-(\pa_2\psi)^\dg(\pa_1\psi)=2\i\Im((\pa_1\psi)^\dg(\pa_2\psi)).
\end{equation}
Now, comparing \eqref{represent1} and \eqref{represent2}, we obtain
\[2\hbar\Im((\pa_1\psi)^\dg\pa_2\psi) = \frac{\hbar}{\i}\tr(P[\pa_1P,\pa_2P])=\frac{\hbar}{2}n\cdot(\pa_1n\times\pa_2 n),\]
and consequently,
\begin{equation}\label{represent_v}
	v=\nabla^\perp\Delta_D^{-1}(2\hbar\Im((\pa_1\psi)^\dg\pa_2\psi))=\nabla^\perp \Delta^{-1}_D\left(\frac{\hbar}{2}n\cdot(\pa_1n\times\pa_2n)\right).
\end{equation}
This completes the derivation of \eqref{eq:ISF_director}.

\subsection{Derivation of \eqref{eq:ISF_velocity}}
Now, we identify the equation that $v$ satisfies. Taking time derivative on $v_j=\hbar\Im(\psi^\dg\pa_j\psi)$, we get
\begin{align*}
	\pa_t v_j &= \hbar\pa_t\Im(\psi^\dg\pa_j\psi)= \hbar \Im(\pa_t\psi^\dg\pa_j\psi)+\hbar\Im(\psi^\dg\pa_j\pa_t\psi)\\
	&=\hbar \Im(\pa_t\psi^\dg\pa_j\psi)+\hbar\pa_j\Im(\psi^\dg\pa_t\psi)-\hbar\Im(\pa_j\psi^\dg\pa_t\psi)\\
	&=\hbar\pa_j\Im(\psi^\dg\pa_t\psi) + 2\hbar\Im(\pa_t\psi^\dg\pa_j\psi).
\end{align*}
The first term can be computed by using \eqref{eq:ISF} and $|\psi|^2=1$ as
\begin{align*}
	\Im(\psi^\dg\pa_t\psi) &=\Im\left(\psi^\dg\left(\frac{\i\hbar}{2}\Delta\psi -\frac{\i}{\hbar}p\psi\right)\right)=\frac{\hbar}{2}\Re(\psi^\dg\Delta\psi)-\frac{p}{\hbar}\\
	&=\frac{\hbar}{2}\pa_j \Re(\psi^\dg\pa_j\psi) - \frac{\hbar}{2}|\nabla\psi|^2-\frac{p}{\hbar}=-\frac{\hbar}{2}|\nabla\psi|^2-\frac{p}{\hbar},
\end{align*}
while the second term can be computed by using \eqref{eq:ISF} again as
\begin{align*}
	\Im(\pa_t\psi^\dg\pa_j\psi) &= \Im\left(\left(\frac{\i\hbar}{2}\Delta\psi-\frac{\i}{\hbar}p\psi\right)^\dg\pa_j\psi\right)=-\frac{\hbar}{2}\Re\left((\Delta \psi)^\dg\pa_j\psi\right)\\
	&=-\frac{\hbar}{2}\pa_k\Re\left((\pa_k\psi)^\dg\pa_j\psi\right)+\frac{\hbar}{2}\Re\left((\pa_k\psi)^\dg\pa_j\pa_k\psi\right)\\
	&=-\frac{\hbar}{2}\pa_k\Re\left((\pa_k\psi)^\dg\pa_j\psi\right)+\frac{\hbar}{4} \pa_j |\nabla\psi|^2.
\end{align*}
Combining these identities, we get
\begin{align}
	\begin{aligned}\label{patvj}
	\pa_tv_j &= \hbar \pa_j\left(-\frac{\hbar}{2}|\nabla\psi|^2-\frac{p}{\hbar}\right)-\hbar^2\pa_k\Re\left((\pa_k\psi)^\dg\pa_j\psi\right)+\frac{\hbar^2}{2}\pa_j|\nabla\psi|^2\\
	&=-\pa_jp -\hbar^2\pa_k\Re\left((\pa_k\psi)^\dg\pa_j\psi\right).
	\end{aligned}
\end{align}
To represent the quantity $\Re((\pa_k\psi)^\dg\pa_j\psi)$ in terms of $v$ and $n$, we again consider $P=\psi\psi^\dg=\frac{1}{2}(I+n\cdot\sigma)$. We use \eqref{Pauli_identity_2} to compute the quantity $(\pa_jP)(\pa_kP)$ as 
\[(\pa_jP)(\pa_kP) = \frac{1}{4}(\pa_jn\cdot\sigma)(\pa_kn\cdot\sigma) = \frac{1}{4}(\pa_j n\cdot\pa_k n)I +\frac{\i}{4}(\pa_jn\times\pa_kn)\cdot\sigma.\]
Taking the trace and using $\tr(\sigma)=0$, we get
\begin{equation}\label{pajPpakP1}
	\tr((\pa_jP)(\pa_kP))=\frac{1}{2}(\pa_jn\cdot\pa_k n).
\end{equation}
On the other hand, since
\[\pa_jP = (\pa_j\psi) \psi^\dg + \psi(\pa_j\psi)^\dg,\]
we get
\begin{align*}
	(\pa_jP)(\pa_kP) &= ((\pa_j\psi)\psi^\dg+\psi(\pa_j\psi)^\dg)((\pa_k\psi)\psi^\dg+\psi(\pa_k\psi)^\dg) \\
	&= (\pa_j\psi)\psi^\dg (\pa_k\psi)\psi^\dg +\psi(\pa_j\psi)^\dg(\pa_k\psi)\psi^\dg + (\pa_j\psi)(\pa_k\psi)^\dg +\psi(\pa_j\psi)^\dg\psi(\pa_k\psi)^\dg.
\end{align*}
Therefore, we again take the trace and use \eqref{psidgpajpsi} to get
\begin{align}
	\begin{aligned}\label{pajpakP2}
	\tr((\pa_jP)(\pa_kP)) &= \psi^\dg(\pa_j\psi)\psi^\dg (\pa_k\psi) +(\pa_j\psi)^\dg(\pa_k\psi) + (\pa_k\psi)^\dg(\pa_j\psi) +(\pa_j\psi)^\dg\psi(\pa_k\psi)^\dg\psi\\
	&=\left(\frac{\i}{\hbar}v_j\right)\left(\frac{\i}{\hbar}v_k\right)+ 2\Re\left((\pa_j\psi)^\dg(\pa_k\psi)\right)+\left(-\frac{\i}{\hbar}v_j\right)\left(-\frac{\i}{\hbar}v_k\right)\\
	&=-\frac{2}{\hbar^2}v_jv_k+2\Re\left((\pa_j\psi)^\dg(\pa_k\psi)\right).
	\end{aligned}
\end{align}
From \eqref{pajPpakP1} and \eqref{pajpakP2}, we obtain
\[\hbar^2\Re\left((\pa_k\psi)^\dg\pa_j\psi\right)=v_jv_k+\frac{\hbar^2}{4}(\pa_j n\cdot\pa_k n),\]
and therefore, returning to \eqref{patvj}, we get
\[\pa_t v_j +\pa_j p = -\pa_k(v_jv_k)-\frac{\hbar^2}{4}\pa_k(\pa_j n\cdot\pa_k n),\]
which is the velocity equation \eqref{eq:ISF_velocity}.

\section{Regularization of \eqref{eq:ISF_director}--\eqref{v:Kn}}\label{sec:3}
\setcounter{equation}{0}

In this section, we introduce the regularized system for \eqref{eq:ISF_director}--\eqref{v:Kn} with parameter $\e$, and present the well-posedness of it for each fixed $\e$. Since the time rescaling $t\mapsto \frac{\hbar}{2}t$ transforms \eqref{eq:ISF_director} into the same system with $\hbar=2$, we assume $\hbar=2$ in Sections \ref{sec:3} and Section \ref{sec:4} without loss of generality. After this normalization, the system \eqref{eq:ISF_director}--\eqref{v:Kn} is simplified as

\begin{equation}\label{eq:ISF_director_h2}
	\pa_t n + (v\cdot\nabla)n =n\times \Delta n,\quad v=K[n]:=\nabla^\perp \Delta^{-1}_D\omega,\quad \omega:=n\cdot(\pa_1n \times\pa_2n).
\end{equation}

Then, we consider the regularized system of \eqref{eq:ISF_director_h2} given as
\begin{equation}\label{eq:ISF_regularize}
	\pa_t n =\e \tau_v(n)+n\times \tau_v(n),
\end{equation}
subject to the initial data $n(0)=n_0$ and the boundary condition $\pa_\nu n=0$ on $\pa\Omega$. Here, $\tau_v(n)$ is a tension field defined as
\[\tau_v(n):=\Delta n +|\nabla n|^2n + n\times(v\cdot\nabla)n.\]
Indeed, we note that whenever $|n|=1$, 
\begin{align*}
	n\times \tau_v(n) &= n\times (\Delta n+|\nabla n|^2n +n\times(v\cdot \nabla)n)=n\times \Delta n +0 + n\times (n\times (v\cdot \nabla)n)\\
	&=n\times \Delta n + n (n\cdot(v\cdot \nabla )n)-(v\cdot\nabla )n (n\cdot n)\\
	&=n\times \Delta n -(v\cdot\nabla) n,
\end{align*}
where the last equality comes from $|n|^2=1$ and 
\begin{equation}\label{nvnablan}
	n\cdot(v\cdot\nabla) n = n_i (v_j \pa_jn_i)=v_j\pa_j|n|^2=0.
\end{equation}
Thus, \eqref{eq:ISF_regularize} reduces to \eqref{eq:ISF_director_h2} when $\e=0$ and therefore it is an appropriate regularized system. The goal of this section is to prove the existence of the solution to \eqref{eq:ISF_regularize}.

\begin{proposition}\label{prop:regularize_exist}
	Let $n_0\in H^3(\Omega;\bbs^2)$ with $\pa_\nu n_0=0$ on $\pa\Omega$. Then, there exists a positive time $T_\e>0$ that only depends on the initial data $n_0$ and $\e$ such that a unique solution $n\in C([0,T_\e];H^3)\cap L^2(0,T_\e;H^4)$ to \eqref{eq:ISF_regularize} satisfying $\pa_tn \in C([0,T_\e];H^1)\cap L^2(0,T_\e;H^2)$ exists. 
\end{proposition}

\begin{remark}
	As we clarified in the introduction, we follow a similar strategy used in \cite{CWW26} to prove Proposition \ref{prop:regularize_exist}, after careful modification, since we consider the $n$-dependent velocity $v=K[n]$, not the prescribed incompressible velocity. 
\end{remark}

\subsection{Preliminaries}
Before we start to prove Proposition \ref{prop:regularize_exist}, we first present several properties of the tension field $\tau_v(n)$ and the velocity map $v=K[n]$.

\begin{lemma}\label{lem:tau}
	Suppose $n\in C([0,T];H^3)$ with $\pa_tn\in C([0,T];H^1)$ satisfies \eqref{eq:ISF_regularize} with a given velocity field $v\in L^\infty(0,T;L^\infty(\Omega))$, together with the boundary condition $\pa_\nu n=0$ on $\pa\Omega$ and $|n(0)|=1$. Let $\tau:=\tau_v(n)$ and $w:=\pa_t n$. Then,
	\begin{enumerate}
		\item $|n(t,x)|=1$ for all $t\in [0,T]$ and $x\in\Omega$.
		\item $\tau\cdot n=0$ and $n\times \tau =n\times\Delta n-(v\cdot \nabla) n$.
		\item The following identities hold:
		\[\tau = \frac{\e w-n\times w}{1+\e^2},\quad w\times \tau = -\frac{|w|^2}{1+\e^2}n,\quad |\tau|^2=\frac{|w|^2}{1+\e^2}.\]
	\end{enumerate}
\end{lemma}

\begin{proof}
	\noindent (1) To show that $|n(t,x)|=1$, we define $m(t,x):=|n(t,x)|^2-1$. Since $n\cdot (n\times \tau)=0$, we have
	\begin{align*}
		\pa_tm &= 2n\cdot\pa_t n =2\e n\cdot\tau=2\e( n\cdot\Delta n +|\nabla n|^2|n|^2)\\
		&=\e(\Delta |n|^2-2|\nabla n|^2+2|\nabla n|^2(1+m))=\e(\Delta m +2|\nabla n|^2m),
	\end{align*}
	where we used the identity $n\cdot\Delta n =\frac{1}{2}\Delta |n|^2-|\nabla n|^2$. Thus, we have
	\[\frac{1}{2}\frac{d}{dt}\|m\|_{L^2}^2 = \int_\Omega m\,\pa_t m\,dx=-\e\|\nabla m\|_{L^2}^2+2\e\int_{\Omega}|\nabla n|^2m^2\,dx\le 2\e \|\nabla n\|_{L^\infty}^2\|m\|_{L^2}^2,\]
	where we used $\pa_\nu m=2n\cdot\pa_\nu n=0$ on $\pa\Omega$ in the integration by parts. Since $\|\nabla n\|_{L^\infty}\le C\|n\|_{H^3}<+\infty$ and $m(0,x)=0$ for all $x\in\Omega$, the Gr\"onwall inequality implies $m(t,x)=0$ for all $t\in[0,T]$ and $x\in\Omega$.\\
	
	\noindent (2) It directly follows from the definition of $\tau$ and $|n|=1$ that
	\[\tau\cdot n = n\cdot \Delta n +|\nabla n|^2=\frac{1}{2}\Delta |n|^2=0.\]
	Furthermore, using the vector identity $a\times(b\times c)= (a\cdot c)b-(a\cdot b)c$, we get
	\begin{align*}
		n\times \tau &= n\times(\Delta n +|\nabla n|^2 n +n\times (v\cdot \nabla)n)\\
		&=n\times \Delta n + 0 + (n\cdot(v\cdot \nabla) n) n - (v\cdot\nabla)n\\
		&=n\times \Delta n -(v\cdot \nabla)n,
	\end{align*}
	where we used \eqref{nvnablan} in the last equality.\\
	
	\noindent (3) Since $w=\pa_t n =\e\tau+n\times \tau$, we use (1) and (2) to get 
	\[n\times w=\e n\times \tau + n\times(n\times \tau)=\e n\times \tau +(n\cdot \tau) n -|n|^2\tau=\e n\times \tau -\tau,\]
	which implies
	\[\e w-n\times w= (\e^2\tau +\e n\times \tau) - (\e n\times \tau -\tau)=(1+\e^2)\tau.\]
	This proves the first identity. Using this, we directly obtain 
	\[w\times \tau =\frac{w\times(\e w-n\times w)}{1+\e^2} = -\frac{|w|^2n-(w\cdot n)w}{1+\e^2}=-\frac{|w|^2n}{1+\e^2}.\]
	Here, we used $0=\pa_t|n|^2 = 2n\cdot\pa_t n= 2n\cdot w$ in the last equality. Finally, we again use the first identity to observe that
	\[|\tau|^2 = \frac{\e^2|w|^2 -2\e w\cdot(n\times w)+|n\times w|^2}{(1+\e^2)^2}=\frac{\e^2|w|^2-0+|w|^2}{(1+\e^2)^2}=\frac{|w|^2}{1+\e^2},\]
	where we used $|n\times w|^2 = |n|^2|w|^2-(n\cdot w)^2=|w|^2$.
\end{proof}

\begin{lemma}\label{lem:velocity_map}
	Let $n\in H^3(\Omega;\bbs^2)$ and $v:=K[n]:=\nabla^\perp \Delta_D^{-1}(n\cdot(\pa_1n\times\pa_2 n))$. Then,
	\begin{enumerate}
		\item $\dv\,v=0$ in $\Omega$ and $v\cdot\nu=0$ on $\pa\Omega$.
		\item For any $p\in(1,\infty)$,
		\begin{align*}
			&\|v\|_{W^{1,p}}\le C\|n\|_{H^2}^2,\quad \|v\|_{L^\infty}\le C\|n\|_{H^2}^2,\\
			&\|\nabla v\|_{L^\infty}\le C\|n\|_{H^2}\|n\|_{H^3},\quad \|v\|_{H^2}\le C\|n\|_{H^2}\|n\|_{H^3}.
		\end{align*}
		\item Suppose $n,\tilde{n}\in H^3(\Omega;\bbs^2)$, with $\|n\|_{H^3},\|\tilde{n}\|_{H^3}\le M$. Let $v:=K[n]$ and $\tilde{v}:=K[\tilde{n}]$. Then,
		\[\|v-\tilde{v}\|_{L^2}\le C(M)\|n-\tilde{n}\|_{L^2},\quad \|v-\tilde{v}\|_{H^1}\le C(M)\|n-\tilde{n}\|_{H^1}.\]
	\end{enumerate}
\end{lemma}

\begin{proof}
	\noindent (1) Since, $\dv\,\nabla^\perp=(\pa_1,\pa_2)\cdot(-\pa_2,\pa_1)\equiv0$,
	\[\dv\,v=\dv\,K[n] = \dv\,\nabla^\perp \Delta^{-1}_D (n\cdot(\pa_1n\times\pa_2n))=0.\]
	Furthermore, if we denote $\phi:=\Delta_D^{-1}(n\cdot(\pa_1n\times\pa_2n))$, then $\phi\equiv0$ on the boundary $\pa\Omega$. Therefore, on $\pa\Omega$,
	\[v\cdot\nu = \nabla^\perp \phi \cdot\nu=-\pa_2\phi \nu_1 +\pa_1\phi\nu_2 = -\pa_{\nu^\perp}\phi= 0,\]
	where $\pa_{\nu^\perp}$ denotes the tangential derivative of $\phi$ on the boundary.\\
	
	\noindent (2) Again, let $\phi:=\Delta_D^{-1}(n\cdot(\pa_1n\times\pa_2n))$. Then, it follows from the standard elliptic estimate $\|\Delta_D^{-1}f\|_{W^{2,p}}\le C\|f\|_{L^p}$ for $p\in(1,\infty)$ that
	\[\|v\|_{W^{1,p}}\le \|\phi\|_{W^{2,p}}\le C\|n\cdot(\pa_1n\times\pa_2 n)\|_{L^p}\le C\|\nabla n\|_{L^{2p}}^2,\]
	where we used the fact $|n|=1$. Since $p\in (1,\infty)$, the Sobolev embedding yields $\|\nabla n\|_{L^{2p}}\le C\|\nabla n\|_{H^1}\le C\|n\|_{H^2}$, which proves the first estimate. The second inequality follows from the first estimate and Sobolev embedding $W^{1,p}\subset L^\infty$ for $p>2$. For the third inequality, we first note that for $p>2$,
	\[\|\nabla v\|_{L^\infty}\le \|\phi\|_{W^{2,\infty}}\le C\|\phi\|_{W^{3,p}}\le C\|n\cdot(\pa_1n\times\pa_2n)\|_{W^{1,p}}.\]
	Since we already obtained for $p>1$ that
	\[\|n\cdot(\pa_1n\times\pa_2 n)\|_{L^p}\le \|\nabla n\|_{L^{2p}}^2\le C\|n\|_{H^2}^2,\]
	we focus on the estimate on the derivative. Note that for $j=1,2$,
	\[\pa_j\left(n\cdot(\pa_1n\times\pa_2 n)\right)=n\cdot(\pa_j\pa_1n\times\pa_2 n)+n\cdot(\pa_1n\times\pa_j\pa_2n).\]
	Therefore, using $|n|=1$, we get for all $p>1$ that 
	\[\|\nabla(n\cdot(\pa_1n\times\pa_2n))\|_{L^p}\le C\|\nabla n\|_{L^{2p}}\|\nabla^2n\|_{L^{2p}}\le C\|n\|_{H^2}\|n\|_{H^3}.\]
	This completes the proof of the third estimate. Finally, the same computation with $p=2$ and the elliptic estimate $\|\Delta_D^{-1}f\|_{H^3}\le C\|f\|_{H^1}$ give
	\[\|v\|_{H^2}\le \|\phi\|_{H^3}\le C\|n\cdot(\pa_1n\times\pa_2 n)\|_{H^1}\le C\|n\|_{H^2}\|n\|_{H^3},\]
	which is the last estimate.\\
	
	\noindent (3) Denote $\omega:=n\cdot (\pa_1n\times\pa_2n)$ and $\tilde{\omega}:=\tilde{n}\cdot(\pa_1\tilde{n}\times\pa_2\tilde{n})$. We first estimate the difference $\omega-\tilde{\omega}$. Define $\delta n:=n-\tilde{n}$. Then,
	\begin{align*}
		\omega-\tilde{\omega} = \delta n \cdot(\pa_1n\times \pa_2n)+\tilde{n} \cdot(\pa_1\delta n\times \pa_2 n)+\tilde{n}\cdot(\pa_1\tilde{n}\times\pa_2 \delta n).
	\end{align*}
	Therefore,
	\[\|\omega-\tilde{\omega}\|_{L^2}\le \|\delta n\|_{L^2}\|\nabla n\|_{L^\infty}^2+\|\nabla \delta n\|_{L^2}\|\nabla n\|_{L^\infty}+\|\nabla \tilde{n}\|_{L^\infty}\|\nabla \delta n\|\le C(M)\|\delta n\|_{H^1}.\]
	This gives
	\[\|v-\tilde{v}\|_{H^1}\le \|\omega-\tilde{\omega}\|_{L^2}\le C(M)\|n-\tilde{n}\|_{H^1}.\]
	To get the $L^2$-estimate, we first note that $\|v-\tilde{v}\|_{L^2}\le \|\nabla \Delta_D^{-1}(\omega-\tilde{\omega})\|_{L^2}\le \|\omega-\tilde{\omega}\|_{H^{-1}}$. Hence, we estimate the $H^{-1}$-norm of $\omega-\tilde{\omega}$ as follows. For any test function $\phi\in H^1_0$, 
	\begin{align*}
		\int_\Omega (\omega-\tilde{\omega})\phi\,dx&=\int_\Omega \delta n\cdot(\pa_1n\times\pa_2n)\phi\,dx+\int_\Omega \tilde{n}\cdot(\pa_1 \delta n\times\pa_2n)\phi\,dx\\
		&\quad +\int_\Omega \tilde{n}\cdot(\pa_1\tilde{n}\times\pa_2 \delta n)\phi\,dx.
	\end{align*}
	The first term is simply bounded by $\|\delta n\|_{L^2}\|\nabla n\|_{L^\infty}^2\|\phi\|_{L^2}$. For the second term, we use integration by parts and $|\tilde{n}|=1$ to get
	\begin{align*}
		\int_\Omega \tilde{n}\cdot(\pa_1\delta n\times \pa_2 n)\phi\,dx &= -\int_\Omega \tilde{n}\cdot(\delta n\times\pa_2n)\pa_1\phi\,dx -\int_\Omega \pa_1\tilde{n}\cdot(\delta n\times\pa_2n)\phi\,dx\\
		&\quad -\int_\Omega \tilde{n}\cdot(\delta n\times\pa_1\pa_2n)\phi\,dx\\
		&\le \|\delta n\|_{L^2}\|\nabla n\|_{L^\infty}\|\nabla\phi\|_{L^2}+\|\nabla \tilde{n}\|_{L^\infty}\|\delta n\|_{L^2}\|\nabla n\|_{L^\infty}\|\phi\|_{L^2}\\
		&\quad + \|\delta n\|_{L^2}\|\nabla^2n\|_{L^4}\|\phi\|_{L^4}\\
		&\le C(M)\|\delta n\|_{L^2} \|\phi\|_{H^1}.
	\end{align*}
	The third term can be estimated by using the same argument as the second term. Thus, we get
	\[\int_\Omega(\omega-\tilde{\omega})\phi\,dx \le C(M)\|\delta n\|_{L^2}\|\phi\|_{H^1},\]
	which implies $\|\omega-\tilde{\omega}\|_{H^{-1}}\le C(M)\|n-\tilde{n}\|_{L^2}$. This proves the desired estimate.
\end{proof}

\subsection{Approximated system}
Now, we begin the proof of Proposition \ref{prop:regularize_exist}. For each $k\in \bbn\cup\{0\}$, we consider a sequence of approximated system of \eqref{eq:ISF_regularize} given as

\begin{equation}\label{eq:iteration}
	\pa_t n^{k+1} = \e\tau_{v^{k}}(n^{k+1})+n^{k+1}\times \tau_{v^{k}}(n^{k+1}),
\end{equation}
subject to initial data $n^{k+1}(0,x)=n_0(x)$ and the boundary condition $\pa_\nu n^{k+1}=0$ on $\pa\Omega$. Here $v^{k}:=K[n^{k}]$ and $n^{0}:=n_0$. To obtain the uniform bound on the sequence $\{n^{k}\}_{k=0}^\infty$, we first consider the equation with frozen velocity:
\begin{equation}\label{eq:frozen}
	\pa_t n = \e\tau_{\bar{v}}(n)+n\times \tau_{\bar{v}}(n),
\end{equation}
where $\bar{v}$ is a given vector field. Before we deal with \eqref{eq:frozen}, we define the solution class $\mathcal{C}_T$ as a set of functions $n:[0,T]\times\Omega\to \bbs^2$ satisfying $\pa_\nu n=0$ on $\pa\Omega$ and
\[n\in C([0,T];H^3)\cap L^2(0,T;H^4),\quad \pa_tn\in C([0,T];H^1)\cap L^2(0,T;H^2).\]

Furthermore, we define the velocity field class $\mathcal{V}_T$ as the set of vector fields $v:[0,T]\times\Omega\to\bbr^2$ satisfying $\dv\,v=0$ in $\Omega$,
$v\cdot\nu=0$ on $\pa\Omega$, and
\[v\in L^\infty(0,T;H^2\cap W^{1,\infty})\cap C([0,T];H^1),\quad
\pa_tv\in L^2(0,T;H^1).\]

The following theorem states that, for a given vector field $\bar{v}\in\mathcal{V}_T$ and initial data $n_0\in H^3(\Omega;\bbs^2)$, there exists a positive time $T_{\e,\bar{v}}\le T$, that only depends on $n_0$, $\bar{v}$ and $\e$ such that \eqref{eq:frozen} has a solution $n$ up to time $T_\e$.

\begin{theorem}\label{thm:frozen_velocity}
	Let $n_0\in H^3(\Omega;\bbs^2)$ with $\pa_\nu n_0=0$ on $\pa\Omega$, and let
	$\bar{v}\in\mathcal{V}_T$. Then, there exists $T_{\e,\bar{v}}\in(0,T]$, that only depends on
	$\|n_0\|_{H^2}$, $\|\bar{v}\|_{L^\infty W^{1,3}}$, and $\e$, such that
	\eqref{eq:frozen} has a unique solution $n\in\mathcal{C}_{T_{\e,\bar{v}}}$.
\end{theorem}

\begin{remark}
	Theorem \ref{thm:frozen_velocity} is proved in \cite[Theorem 4.1]{CWW26} for a
	three-dimensional domain, under the weaker assumptions
	$\bar v\in L^\infty(0,T;W^{1,3})\cap C([0,T];H^1)$, $\pa_t\bar v\in L^2(0,T;H^1)$,
	$\dv\,\bar v=0$ in $\Omega$ and $\bar v\cdot\nu=0$ on $\pa\Omega$, which are implied
	by $\bar v\in\mathcal V_T$ since $\Omega$ is bounded. Therefore its proof applies to our
	two-dimensional domain without change. 
\end{remark}

When we apply Theorem \ref{thm:frozen_velocity} to construct the sequence of approximated solutions $\{n^{k}\}_{k=0}^\infty$ to \eqref{eq:iteration}, we need to pay attention on the regularity and the norm of $v^{k}$. Precisely, we first need to verify that the velocity $v^{k}$ belongs to $\mathcal{V}_T$, so that we use it to construct the solution $n^{k+1}$ at the next step. Furthermore, the existence time $T_{\e,v^k}$ on which the solution $n^{k+1}$ is defined depends on the norm $\|v^{k}\|_{L^\infty W^{1,3}}$. Therefore, without uniform-in-$k$ estimate on $\|v^{k}\|_{L^\infty W^{1,3}}$, the existence time of $n^k$ may shrink to 0 as $k\to\infty$, preventing to construct the solution $n$ of \eqref{eq:ISF_regularize} as a limit of $n^k$ on the positive time interval. Therefore, we need to show that $\|v^{k}\|_{L^\infty W^{1,3}}$ is uniformly bounded so that the existence time is uniformly bounded from below. We prove these two things for the solution $n$ to the frozen velocity equation \eqref{eq:frozen} in the following subsections.

\subsection{Regularity of the velocity field}

We first show that, the velocity field $v=K[n]$ generated from the solution $n$ to \eqref{eq:frozen} has enough regularity to apply Theorem \ref{thm:frozen_velocity} again with that velocity $v$.

\begin{lemma}\label{lem:v}
	For $n\in\mathcal{C}_T$, the velocity field $v=K[n]$ belongs to $\mathcal{V}_T$.
\end{lemma}

\begin{proof}
	It directly follows from Lemma \ref{lem:velocity_map} (1) that $v$ satisfies the divergence-free condition and slip boundary condition. Moreover, by Lemma \ref{lem:velocity_map} (2) and the assumption $n\in\mathcal{C}_T$, we get
	\begin{align*}
		&\|v\|_{L^\infty L^\infty}\le C\|n\|_{L^\infty H^2}^2<+\infty,\\
		&\|v\|_{L^\infty H^2} \le C\|n\|_{L^\infty H^2}\|n\|_{L^\infty H^3}<+\infty,\\
		&\|\nabla v\|_{L^\infty L^\infty}\le C\|n\|_{L^\infty H^2}\|n\|_{L^\infty H^3}<+\infty.
	\end{align*}
	Next, to show $\pa_t v\in L^2(0,T;H^1)$, we first define $\omega:=n\cdot(\pa_1n\times\pa_2 n)$. Then, 
	\[\|\pa_t v\|_{H^1}\le \|\Delta^{-1}_D\pa_t\omega\|_{H^2} \le C\|\pa_t \omega\|_{L^2}.\]
	On the other hand, since $\pa_tn$, $\pa_1n$, and $\pa_2n$ are vectors in $\bbr^3$ that are perpendicular to $n$, they are linearly dependent and therefore, $\pa_tn\cdot(\pa_1n\times\pa_2n)=0$. Hence,
	\begin{align*}
		\pa_t \omega &= \pa_tn \cdot(\pa_1n\times\pa_2n)+n\cdot(\pa_1\pa_tn\times\pa_2 n)+n\cdot(\pa_1n\times\pa_2\pa_tn)\\
		&=n\cdot(\pa_1\pa_tn\times\pa_2 n+\pa_1n\times\pa_2\pa_tn).
	\end{align*} 
	Since $|n|=1$, this shows 
	\[\|\pa_t\omega\|_{L^2}\le 2\|\nabla\pa_tn\|_{L^2}\|\nabla n\|_{L^\infty}\le C\|n\|_{H^3}\|\pa_tn\|_{H^1},\]
	and consequently, thanks to the assumption $n\in\mathcal{C}_T$,
	\[\|\pa_t v\|_{L^2H^1}\le C\|\pa_t\omega\|_{L^2L^2}\le C\|n\|_{L^2H^3}\|\pa_tn\|_{L^\infty H^1}<+\infty.\]
	Finally, to show that $v\in C([0,T];H^1)$. we first show that $\omega\in C([0,T];L^2)$. Indeed, for any $t,s\in [0,T]$,
	\begin{align*}
		\omega(t)-\omega(s)&=(n(t)-n(s))\cdot(\pa_1n(t)\times\pa_2n(t))\\
		&\quad +n(s)\cdot(\pa_1(n(t)-n(s))\times \pa_2n(t))\\
		&\quad +n(s)\cdot(\pa_1n(s)\times\pa_2(n(t)-n(s))).
	\end{align*}
	Therefore,
	\begin{align*}
		\|\omega(t)-\omega(s)\|_{L^2}&\le \|n(t)-n(s)\|_{L^2}\|\nabla n(t)\|_{L^\infty}^2+\|\nabla(n(t)-n(s))\|_{L^2}(\|\nabla n(t)\|_{L^\infty}+\|\nabla n(s)\|_{L^\infty})\\
		&\le C(\|n(t)\|_{H^3}^2+\|n(t)\|_{H^3}+\|n(s)\|_{H^3})\|n(t)-n(s)\|_{H^1},
	\end{align*}
	which implies that $\omega\in C([0,T];L^2)$ since $n\in C([0,T];H^3)$. From
	\[\|v(t)-v(s)\|_{H^1}\le\|\Delta_D^{-1}(\omega(t)-\omega(s))\|_{H^2}\le C\|\omega(t)-\omega(s)\|_{L^2},\]
	we conclude $v\in C([0,T];H^1)$. 
\end{proof}

\subsection{Invariance of the norm} 

Next, we prove that there exists a positive time $T_\e$ such that, if $\bar{n}\in\mathcal{C}_{T_\e}$ is bounded in a particular norm, then the solution $n$ to \eqref{eq:frozen} belongs to the same solution class $\mathcal{C}_{T_\e}$ satisfying the same norm estimate. Before we state this property, we first present several lemmas that are required in the proof of norm invariance.

\begin{lemma}\label{lem:norm_equiv}
	Suppose $n\in\mathcal{C}_T$ solves \eqref{eq:frozen} with $\bar{v}\in\mathcal{V}_T$. Then, for $\Gamma:[0,T]\to\bbr$ defined as 
	\begin{equation}\label{eq:Gamma}
		\Gamma(t):=1+\|n(t)\|_{H^2}+\|\pa_tn(t)\|_{H^1}+\|\bar{v}(t)\|_{L^\infty}+\|\bar{v}(t)\|_{W^{1,4}},
	\end{equation} 
	there exists a positive constant $C$ such that 
	\[\|\nabla\Delta n(t)\|_{L^2}\le C \Gamma(t)^3,\quad \|n(t)\|_{H^3}\le C \Gamma(t)^3.\]
\end{lemma}

\begin{proof}
	It follows from the definition of $\tau=\tau_{\bar{v}}(n)$ that 
	\[\nabla\Delta n = \nabla\tau -\nabla(|\nabla n|^2n)-\nabla(n\times(\bar{v}\cdot\nabla)n).\]
	We estimate the $L^2$-norms of the three terms on the right-hand side of the above equation. First, using Lemma \ref{lem:tau} (3) to get
	\[\pa_j\tau = \frac{\e\pa_jw -\pa_j n\times w -n\times \pa_jw}{1+\e^2},\quad j=1,2,\]
	where $w = \pa_t n$. Then, we have
	\begin{align*}
		\|\nabla\tau\|_{L^2} &\le C\|\nabla\pa_t n\|_{L^2} +\|\nabla n\|_{L^4}\|\pa_tn\|_{L^4}\le C\Gamma+C\Gamma^2.
	\end{align*}
	Next, observing that $\pa_j(|\nabla n|^2n)=2(\pa_j\pa_kn\cdot\pa_kn)n+|\nabla n|^2\pa_jn$, we use the Sobolev interpolation to derive
	\begin{align*}
		\|\nabla(|\nabla n|^2n)\|_{L^2}&\le 2\|\nabla^2n\|_{L^4}\|\nabla n\|_{L^4}+\|\nabla n\|_{L^6}^3\\
		&\le C\|\nabla^2n\|_{L^2}^{1/2}\|\nabla^2n\|_{H^1}^{1/2}\|\nabla n\|_{L^4}+\|\nabla n\|_{L^6}^3\\
		&\le C\Gamma^{3/2}\|n\|_{H^3}^{1/2}+C\Gamma^3.
	\end{align*}
	Finally, we expand the third term as
	\[\pa_j(n\times \bar{v}_i\pa_i n)=\pa_j n\times \bar{v}_i\pa_i n+n\times \pa_j\bar{v}_i\pa_i n +n\times \bar{v}_i\pa_j\pa_in\]
	to estimate it as
	\[\|\nabla(n\times (\bar{v}\cdot\nabla)n)\|_{L^2}\le \|\bar{v}\|_{L^\infty}\|\nabla n\|_{L^4}^2+\|\nabla \bar{v}\|_{L^4}\|\nabla n\|_{L^4}+\|\bar{v}\|_{L^\infty}\|\nabla^2n\|_{L^2}\le C\Gamma^3.\]
	Combining the above estimates, we get
	\begin{equation}\label{est:nableDelta-0}
		\|\nabla \Delta n\|_{L^2}\le C\Gamma^3+C\Gamma^{3/2}\|n\|_{H^3}^{1/2}.
	\end{equation}
	On the other hand, since $|n|=1$ and $\pa_\nu n=0$ on $\pa\Omega$, we have the following elliptic estimate (see, e.g., \cite[Chapter 2, Theorem 5.1]{LM72}):
	\begin{equation}\label{est:elliptic}
		\|n\|_{H^3}\le C(\|n\|_{L^2}+\|\Delta n\|_{H^1})\le C(\Gamma+\|\nabla\Delta n\|_{L^2}).
	\end{equation}
	Therefore, returning to \eqref{est:nableDelta-0}, we obtain
	\[\|\nabla \Delta n\|_{L^2}\le C\Gamma^3+C\Gamma^{3/2}(\Gamma^{1/2}+\|\nabla\Delta n\|^{1/2}_{L^2})\le C\Gamma^3 +\frac{1}{2}\|\nabla\Delta n\|_{L^2}.\]
	This yields the desired bound of $\|\nabla\Delta n\|_{L^2}$. Consequently, \eqref{est:elliptic} implies the required estimate on $\|n\|_{H^3}$.
\end{proof}

\begin{lemma}\label{lem:w}
	Suppose $n\in\mathcal{C}_T$ solves \eqref{eq:frozen} with $\bar{v}\in\mathcal{V}_T$. Let $w:=\pa_t n$. Then, $w$ satisfies
	\begin{align}
		\begin{aligned}\label{eq:w-0}
			\pa_tw &= \e\Delta w+n\times \Delta w -(\bar{v}\cdot\nabla)w-(\pa_t\bar{v}\cdot\nabla)n-(w\cdot(\bar{v}\cdot\nabla)n)n\\
			&\quad -\frac{|w|^2}{1+\e^2}n+|\nabla n|^2(n\times w)+\e R,
		\end{aligned}
	\end{align}
	and $\pa_\nu w=0$, where 
	\begin{equation}\label{def:R}
		R:=|\nabla n|^2w+2(\nabla w:\nabla n)n+w\times (\bar{v}\cdot \nabla)n +n\times (\bar{v}\cdot\nabla)w+n\times (\pa_t\bar{v}\cdot\nabla)n.
	\end{equation}
\end{lemma}

\begin{proof}
	Again, we use simple notation $\tau=\tau_{\bar{v}}$. Since $n\in\mathcal{C}_T$, and in particular $\pa_\nu n=0$ for all $t>0$, we obtain $\pa_\nu w=0$ for all $t$. Hence, we focus on verifying the governing equation of $w$. First, we differentiate \eqref{eq:frozen} by $t$ to get
	\begin{equation}\label{eq:w}
		\pa_t w = \e\pa_t \tau +w\times\tau +n\times\pa_t\tau.
	\end{equation}
	It follows from the definition of $\tau$ that
	\[\pa_t \tau = \Delta w + \pa_t(|\nabla n|^2n)+\pa_t(n\times(\bar{v}\cdot\nabla)n).\]
	Since $\pa_t |\nabla n|^2 = 2\pa_j n\pa_j w=2\nabla n:\nabla w$, and 
	\[\pa_t (n\times (\bar{v}\cdot\nabla)n)=w\times (\bar{v}\cdot\nabla n)+n\times(\pa_t\bar{v}\cdot\nabla)n+n\times(\bar{v}\cdot\nabla)w,\]
	we rewrite $\pa_t \tau$ as
	\[\pa_t \tau = \Delta w+2(\nabla n:\nabla w)n +|\nabla n|^2 w+w\times (\bar{v}\cdot\nabla)n+n\times(\pa_t\bar{v}\cdot\nabla)n+n\times(\bar{v}\cdot\nabla)w.\]
	This yields
	\begin{align*}
		n\times\pa_t \tau &=n\times \Delta w + 0 + |\nabla n|^2 (n\times w) + (n\cdot (\bar{v}\cdot\nabla)n)w-(n\cdot w)(\bar{v}\cdot\nabla)n\\
		&\quad + (n\cdot (\pa_t\bar{v}\cdot\nabla)n)n - (\pa_t\bar{v}\cdot\nabla)n + (n\cdot (\bar{v}\cdot\nabla)w)n - (\bar{v}\cdot\nabla)w\\
		&=n\times\Delta w + |\nabla n|^2(n\times w) + 0 - 0 + 0 -(\pa_t\bar{v}\cdot\nabla)n +(n\cdot(\bar{v}\cdot\nabla)w)n-(\bar{v}\cdot\nabla)w\\
		&=n\times\Delta w + |\nabla n|^2(n\times w) -(\pa_t\bar{v}\cdot\nabla)n -(w\cdot(\bar{v}\cdot\nabla)n)n-(\bar{v}\cdot\nabla)w,
	\end{align*}
	where we used 
	\[n\cdot(\bar{v}\cdot\nabla)w=n\cdot \bar{v}_j(\pa_j w)=\bar{v}_j(\pa_j (n\cdot w)-\pa_j n\cdot w)=-\bar{v}_j\pa_j n \cdot w = -w\cdot(\bar{v}\cdot\nabla)n.\]
	We substitute the above representation into \eqref{eq:w}, and use Lemma \ref{lem:tau} to get the desired equation for $w$.
\end{proof}

Based on the previous lemmas, we show that there exists a positive time $T_\e$, so that if $\bar{n}\in \mathcal{C}_{T_\e}$ is bounded in some appropriate norm, the solution $n\in\mathcal{C}_{T_\e}$ to \eqref{eq:frozen} satisfies the same boundedness.

\begin{proposition}\label{prop:invariance}
	There exist positive constants $Y_*$, $M_*$, and $T_\e$ such that the following holds. Suppose $\bar{n}\in \mathcal{C}_{T_\e}$ satisfies
	\[\sup_{0\le t\le T_\e} \|\Delta \bar{n}(t)\|_{L^2}^2+\|\pa_t \bar{n}(t)\|_{H^1}^2\le Y_*,\quad \sup_{0\le t\le T_\e}\|\bar{n}(t)\|_{H^3}\le M_*,\]
	and $\bar{n}(0)=n_0$. Let $\bar{v}=K[\bar{n}]\in\mathcal{V}_{T_\e}$ be a velocity field generated by $\bar{n}$. Then the solution $n$ to \eqref{eq:frozen} belongs to the same class $\mathcal{C}_{T_\e}$ and satisfies the same boundedness condition:
	\begin{equation}\label{est:bound}
		\sup_{0\le t\le T_\e} \|\Delta n(t)\|_{L^2}^2+\|\pa_t n(t)\|_{H^1}^2\le Y_*,\quad \sup_{0\le t\le T_\e}\|n(t)\|_{H^3}\le M_*.
	\end{equation}
\end{proposition}

\begin{proof}
	We first define 
	\[\bar{y}(t):=\|\Delta \bar{n}(t)\|_{L^2}^2+\|\pa_t \bar{n}(t)\|_{H^1}^2,\quad y(t):=\|\Delta n(t)\|_{L^2}^2+\|\pa_t n(t)\|_{H^1}^2.\]
	Let $v_0:=K[n_0]$ and define
	\[w_0:=\e\tau_{v_0}(n_0)+n_0\times\tau_{v_0}(n_0)\in H^1(\Omega),\]
	which depends only on $n_0$ and $\e$. Furthermore, let
	\[y(0):=\|\Delta n_0\|_{L^2}^2+\|w_0\|_{H^1}^2=:y_0\]
	which is determined by $n_0$ and $\e$, and set $Y_*:=2(1+y_0)$. Also, define $M_*:=C(1+\sqrt{Y_*}+C(1+Y_*))^3$, where $C$ denotes the larger of the constants appearing in Lemma \ref{lem:norm_equiv} and in \eqref{est:barv1} below. Since $|\bar{n}|=1$, with $\pa_\nu \bar{n}=0$ on $\pa\Omega$, it follows from the elliptic estimate that 
	\[\|\bar{n}\|_{H^2} \le C(\|\bar{n}\|_{L^2}+\|\Delta \bar{n}\|_{L^2})\le  C(1+\|\Delta \bar{n}\|_{L^2}).\]
	Thus, together with the assumption, we have
	\[\|\bar{n}\|_{H^2}\le C(1+\sqrt{Y_*}),\quad \|\pa_t\bar{n}\|_{H^1}\le \sqrt{Y_*}.\]
	Therefore, following the proof of Lemma \ref{lem:v}, we get the bounds on $\bar{v}$ as
	\begin{equation}\label{est:barv1}
		\|\bar{v}\|_{L^\infty}+\|\bar{v}\|_{W^{1,4}}\le C\|\bar{n}\|_{H^2}^2\le C(1+Y_*),
	\end{equation}
	\begin{equation}\label{est:barv2}
		\|\pa_t \bar{v}\|_{H^1}\le C\|\bar{n}\|_{H^3}\|\pa_t\bar{n}\|_{H^1}\le CM_*\sqrt{Y_*},
	\end{equation}
	and
	\begin{equation}\label{est:barv3}
		\|\nabla \bar{v}\|_{L^\infty}\le C\|\bar{n}\|_{H^2}\|\bar{n}\|_{H^3}\le C(1+\sqrt{Y_*})M_*.
	\end{equation}
	Since $\|\bar{v}\|_{L^\infty W^{1,3}}\le C\|\bar{v}\|_{L^\infty W^{1,4}}\le C(1+Y_*)$ by \eqref{est:barv1}, Theorem \ref{thm:frozen_velocity} provides the solution $n\in\mathcal{C}_{T'}$ to \eqref{eq:frozen} on an interval $[0,T']$ with $T'>0$ depending only on $\|n_0\|_{H^2}$, $Y_*$ and $\e$. We always take $T_\e\le T'$, so that $n\in\mathcal{C}_{T_\e}$.
	
	To show that the solution $n$ satisfies the same bounds as in \eqref{est:bound}, we estimate the evolution of $y(t)$. Using the definition of $\Gamma$ in \eqref{eq:Gamma}, we get
	\[\frac{d}{dt}\|\Delta n\|_{L^2}^2 = 2\int_\Omega\Delta n\cdot \Delta \pa_t n\,dx=-2\int_\Omega \nabla\Delta n:\nabla \pa_t n\,dx\le 2\|\nabla\Delta n\|_{L^2}\|\nabla \pa_t n\|_{L^2}\le C\Gamma^4.\]
 	Next, we control $\|\pa_tn\|_{H^1}$. For simplicity, we define $w:=\pa_t n$. Then, multiplying \eqref{eq:w-0} with $w$, we get
	\begin{align*}
		\frac{1}{2}\frac{d}{dt}\|w\|_{L^2}^2&=\e\int_\Omega w\cdot\Delta w\,dx+\int_\Omega w\cdot (n\times\Delta w)\,dx-\int_\Omega w\cdot(\bar{v}\cdot\nabla w)\,dx-\int_\Omega w\cdot (\pa_t\bar{v}\cdot\nabla n)\,dx\\
		&\quad -\int_\Omega (w\cdot n)(w\cdot(\bar{v}\cdot\nabla n))\,dx-\frac{1}{1+\e^2}\int_\Omega |w|^2n\cdot w\,dx+\int_\Omega |\nabla n|^2w\cdot(n\times w)\,dx\\
		&\quad +\e\int_\Omega w\cdot R\,dx\\
		&=:\sum_{i=1}^8 I_i.
	\end{align*}
	
	It is straightforward that $I_1 = -\e\|\nabla w\|^2$. Since $w\cdot n = \pa_tn\cdot n =\frac{1}{2}\frac{d}{dt}|n|^2=0$, we get $I_5=I_6=0$. Also, since $w\cdot(n\times w)=0$, we have $I_7=0$. The term $I_3$ also vanishes due to the incompressibility of $\bar{v}$. Thus, we focus on estimating $I_2$, $I_4$, and $I_8$.\\
	
	\noindent $\bullet$ (Estimate of $I_2$): We rewrite $I_2$ as
	\begin{align*}
		I_2&=\int_\Omega w\cdot(n\times \Delta w)\,dx=\int_\Omega \Delta w\cdot(w\times n)\,dx\\
		&=-\int_\Omega\pa_jw\cdot(\pa_j w\times n)\,dx-\int_\Omega\pa_jw \cdot(w\times \pa_jn)\,dx=-\int_\Omega \pa_j w\cdot(w\times \pa_jn)\,dx,
	\end{align*}
	which implies
	\[I_2\le \|\nabla w\|_{L^2}\|w\|_{L^4}\|\nabla n\|_{L^4}\le \Gamma^3.\]
	\noindent $\bullet$ (Estimate of $I_4$): The term $I_4$ can be simply estimated as
	\[I_4\le \|w\|_{L^2}\|\pa_t \bar{v}\|_{L^4}\|\nabla n\|_{L^4}\le C\|\pa_t \bar{v}\|_{L^4}\Gamma^2.\]
	\noindent $\bullet$ (Estimate of $I_8$): Finally, we use \eqref{def:R} and $n\cdot w=0$ to estimate $I_8$ as
	\begin{align*}
		I_8&=\int_\Omega |\nabla n|^2 |w|^2\,dx+ 0+0 +\int_\Omega w\cdot(n\times (\bar{v}\cdot\nabla w))\,dx+\int_\Omega w\cdot n\times (\pa_t \bar{v}\cdot \nabla n)\,dx\\
		&\le \|\nabla n\|_{L^4}^2\|w\|_{L^4}^2+\|w\|_{L^2}\|\bar{v}\|_{L^\infty}\|\nabla w\|_{L^2}+\|w\|_{L^2}\|\pa_t \bar{v}\|_{L^4}\|\nabla n\|_{L^4}\\
		&\le C\Gamma^4+ C\Gamma^3+C\|\pa_t \bar{v}\|_{L^4}\Gamma^2.
	\end{align*}
	Therefore, combining all the estimates on $I_i$ for $i=1,2,\ldots 8$, we get
	\[\frac{1}{2}\frac{d}{dt}\|w\|_{L^2}^2+\e\|\nabla w\|_{L^2}^2\le C(\Gamma^4+\|\pa_t \bar{v}\|_{L^4}\Gamma^2).\]
	Finally, we estimate $\|\nabla w\|_{L^2}$. First, since
	\[\frac{1}{2}\frac{d}{dt}\|\nabla w\|^2_{L^2}=\int_\Omega \nabla w:\nabla\pa_t w\,dx=-\int_\Omega \Delta w\cdot\pa_t w\,dx,\]
	we substitute the equation \eqref{eq:w-0} for $w$ to get
	\begin{align*}
		\frac{1}{2}\frac{d}{dt}\|\nabla w\|_{L^2}^2&=-\e\|\Delta w\|_{L^2}^2+0+\int_\Omega \Delta w\cdot(\bar{v}\cdot\nabla w)\,dx+\int_\Omega \Delta w\cdot (\pa_t\bar{v}\cdot\nabla n)\,dx\\
		&\quad +\int_\Omega(\Delta w\cdot n)(w\cdot(\bar{v}\cdot\nabla n))\,dx+\frac{1}{1+\e^2}\int_\Omega \Delta w\cdot n|w|^2\,dx\\
		&\quad -\int_\Omega |\nabla n|^2\Delta w\cdot(n\times w)\,dx-\e\int_\Omega \Delta w\cdot R\,dx\\
		&=-\e\|\Delta w\|_{L^2}^2+0+\sum_{i=1}^6J_i,
	\end{align*}
	where the second term vanishes due to the incompressibility of $\bar{v}$.\\
	
	\noindent $\bullet$ (Estimate of $J_1$): Again, using incompressibility of $\bar{v}$, we get
	\begin{align*}
		J_1 &= \int_\Omega \bar{v}_i\pa_i w\cdot \pa_{jj}w\,dx=-\int_\Omega \pa_j\bar{v}_i\pa_iw\cdot\pa_j w\,dx-\int_\Omega \bar{v}_i\pa_{ij}w\cdot\pa_jw\,dx\\
		&=-\int_\Omega\pa_j\bar{v}_i\pa_iw\cdot\pa_jw\,dx\le \|\nabla \bar{v}\|_{L^\infty}\|\nabla w\|^2\le \|\nabla \bar{v}\|_{L^\infty}\Gamma^2.
	\end{align*}
	\noindent $\bullet$ (Estimate of $J_2$): We directly estimate $J_2$ as
	\begin{align*}
		J_2&=\int_\Omega \pa_t\bar{v}_i\pa_i n \cdot \pa_{jj}w\,dx=-\int_\Omega\pa_t\pa_j \bar{v}_i\pa_in\cdot\pa_j w\,dx-\int_\Omega\pa_t\bar{v}_i\pa_{ij}n\cdot\pa_jw\,dx\\
		&\le \|\nabla w\|(\|\nabla\pa_tv\|_{L^2}\|\nabla n\|_{L^\infty}+\|\pa_t\bar{v}\|_{L^4}\|\nabla^2n\|_{L^4})\\
		&\le C\Gamma(\|\pa_t \bar{v}\|_{H^1}\Gamma^3+\|\pa_t \bar{v}\|_{H^1}\Gamma^2)\le C\|\pa_t \bar{v}\|_{H^1}\Gamma^4.
	\end{align*}
	\noindent $\bullet$ (Estimate of $J_3$): Using $\pa_j(w\cdot n)=0$, we have 
	\begin{align*}
		J_3 &= -\int_\Omega (\pa_j w\cdot n)\pa_j(w\cdot(\bar{v}\cdot \nabla n))\,dx-\int_\Omega (\pa_j w\cdot\pa_j n)(w\cdot(\bar{v}\cdot \nabla n))\,dx\\
		&= \int_\Omega (w\cdot\pa_jn)\pa_j(w\cdot(\bar{v}\cdot \nabla n))\,dx-\int_\Omega (\pa_j w\cdot\pa_j n)(w\cdot(\bar{v}\cdot \nabla n))\,dx.
	\end{align*}
	By expanding the first term, we estimate it as
	\begin{align*}
		\int_\Omega &(w\cdot\pa_j n)\left(\pa_j w\cdot(\bar{v}\cdot\nabla n)+w\cdot(\pa_j \bar{v}\cdot\nabla n)+w\cdot(\bar{v}\cdot\nabla\pa_j n)\right)\,dx\\
		&\le \|\bar{v}\|_{L^\infty}\|\nabla w\|_{L^2}\|\nabla n\|_{L^8}^2\|w\|_{L^4}+\|\nabla \bar{v}\|_{L^2}\|w\|_{L^8}^2\|\nabla n\|_{L^8}^2+\|\bar{v}\|_{L^\infty}\|w\|_{L^8}^2\|\nabla^2n\|_{L^2}\|\nabla n\|_{L^4}\\
		&\le C\Gamma^5.
	\end{align*}
	On the other hand, the second term can be bounded as
	\[\int_\Omega (\pa_j w\cdot\pa_j n)(w\cdot(\bar{v}\cdot \nabla n))\,dx\le \|\bar{v}\|_{L^\infty}\|w\|_{L^4}\|\nabla n\|_{L^8}^2\|\nabla w\|_{L^2}\le C\Gamma^5.\]
	Thus, we get $J_3\le C\Gamma^5$. \\
	
	\noindent $\bullet$ (Estimate of $J_4$): Using integration by parts and $\pa_j w\cdot n=-w\cdot\pa_j n$, we get
	\begin{align*}
		J_4 &= -\frac{1}{1+\e^2}\int_\Omega \pa_jw \cdot \pa_j(n|w|^2)\,dx=\frac{1}{1+\e^2}\int_\Omega\left(2(w\cdot\pa_jw)(w\cdot\pa_jn)-|w|^2\pa_jw\cdot\pa_jn\right)\,dx\\
		& \le C\|w\|_{L^8}^2\|\nabla w\|_{L^2}\|\nabla n\|_{L^4}\le C\Gamma^4.
	\end{align*}	
	
	\noindent $\bullet$ (Estimate of $J_5$): Applying integration by parts once again and $|n|=1$, we obtain
	\begin{align*}
		J_5&=\int_\Omega \pa_jw \cdot\left(\pa_j(|\nabla n|^2)(n\times w)+|\nabla n|^2(\pa_jn\times w+n\times\pa_jw)\right)\,dx\\
		&=\int_\Omega \pa_jw \cdot\left(\pa_j(|\nabla n|^2)(n\times w)+|\nabla n|^2(\pa_jn\times w)\right)\,dx\\
		&\le C\|\nabla w\|\|\nabla^2 n\|_{L^4}\|\nabla n\|_{L^8}\|w\|_{L^8}+C\|\nabla w\|_{L^2}\|\nabla n\|_{L^{12}}^3\|w\|_{L^4}\le C\Gamma^5.
	\end{align*}
	
	\noindent $\bullet$ (Estimate of $J_6$): Finally, for $J_6$, we first note that $J_6\le \frac{\e}{2}\|\Delta w\|_{L^2}^2+\frac{\e}{2}\|R\|^2_{L^2}$, and
	\begin{align*}
		\|R\|_{L^2}&\le \|\nabla n\|_{L^8}^2\|w\|_{L^4}\\
		&\quad +C\|\nabla n\|_{L^\infty}\|\nabla w\|_{L^2}+\|w\|_{L^4}\|\bar{v}\|_{L^\infty}\|\nabla n\|_{L^4}+\|\bar{v}\|_{L^\infty}\|\nabla w\|_{L^2}+\|\pa_t \bar{v}\|_{L^4}\|\nabla n\|_{L^4}\\
		&\le C(\Gamma^4+\|\pa_t \bar{v}\|_{L^4}\Gamma).
	\end{align*}
	This gives
	\[J_6\le \frac{\e}{2}\|\Delta w\|_{L^2}^2+C(\Gamma^8+\|\pa_t \bar{v}\|_{L^4}^2\Gamma^2).\]
	Combining all the estimates for $J_i$ for $i=1,2,\ldots, 6$, we get
	\begin{align*}
		\frac{1}{2}\frac{d}{dt}\|\nabla w\|_{L^2}^2+\frac{\e}{2}\|\Delta w\|^2_{L^2}&\le C(\|\nabla \bar{v}\|_{L^\infty}\Gamma^2+\|\pa_t \bar{v}\|_{H^1}\Gamma^4+\Gamma^8+\|\pa_t\bar{v}\|_{H^1}^2\Gamma^2)\\
		&\le C((\|\nabla \bar{v}\|_{L^\infty}+\|\pa_t\bar{v}\|_{H^1}^2)\Gamma^2+\Gamma^8).
	\end{align*}
	Therefore, collecting all the estimates for $\|\Delta n\|_{L^2}$ and $\|\pa_t n\|_{H^1}$, we get 
	\begin{equation}\label{est:y}
		\frac{dy}{dt} + 2\e\|\nabla w\|_{L^2}^2+\e\|\Delta w\|_{L^2}^2\le C\Gamma^8+C(\|\nabla \bar{v}\|_{L^\infty}+\|\pa_t\bar{v}\|_{H^1}^2)\Gamma^2.
	\end{equation}
	Finally, since 
	\[\|n\|_{H^2}\le C(1+\|\Delta n\|_{L^2}),\quad \|w\|_{H^1}\le C(\|w\|_{L^2}^2+\|\nabla w\|_{L^2}^2)^{1/2},\]
	we bound $\Gamma$ as
	\[\Gamma\le C(1+\sqrt{y}+\|\bar{v}\|_{L^\infty}+\|\bar{v}\|_{W^{1,4}}),\]
	and therefore, from \eqref{est:y}, we get
	\begin{equation}\label{est:dydt}
		\frac{dy}{dt}\le C(1+y+\|\bar{v}\|^2_{L^\infty}+\|\bar{v}\|^2_{W^{1,4}})^4+C(\|\nabla \bar{v}\|_{L^\infty}+\|\pa_t\bar{v}\|_{H^1}^2)(1+y+\|\bar{v}\|^2_{L^\infty}+\|\bar{v}\|^2_{W^{1,4}}).
	\end{equation}
	Applying \eqref{est:barv1}, \eqref{est:barv2}, and \eqref{est:barv3}, we conclude that there exist constants $a$, $b=C(1+Y_*)$, and $c=C((1+\sqrt{Y_*})M_*+M_*^2Y_*)$ such that the following inequality holds:
	\begin{align*}
		\frac{dy}{dt}\le a(1+y+b^2)^4 +c(1+y+b^2).
	\end{align*}
	Then, $Y(t):=e^{-ct}(1+y(t)+b^2)$ satisfies
	\[\frac{dY}{dt}\le ae^{3ct}Y^4,\]
	from which we obtain 
	\[1+y+b^2\le e^{ct}Y_0\left(1-\frac{a}{c}Y_0^3(e^{3ct}-1)\right)^{-1/3}\quad\mbox{for}\quad t<\frac{1}{3c}\log\left(1+\frac{c}{aY_0^3}\right),\]
	with $Y_0=Y(0)=(1+y_0+b^2)$ and $y_0:=y(0)$. Therefore, if we choose $T_\e$ further small enough so that
	\[e^{cT_\e}\left(1-\frac{a}{c}Y_0^3(e^{3cT_\e}-1)\right)^{-1/3}<1+\frac{2+y_0}{1+y_0+b^2},\]
	then for $0\le t\le T_\e$,
	\begin{align*}
		1+y(t)+b^2&\le Y_0\left(1+\frac{2+y_0}{1+y_0+b^2}\right)=1+y_0+b^2+2+y_0\\
		&=1+b^2+2(1+y_0)=1+b^2+Y_*.
	\end{align*}
	Thus, we get $\sup_{0\le t\le T_\e} y(t)\le Y_*$, which yields the first bound in \eqref{est:bound}. To estimate $H^3$-norm of $n$, we use Lemma \ref{lem:norm_equiv} and \eqref{est:barv1} to get 
	\[\|n\|_{H^3}\le C\Gamma^3\le C(1+\sqrt{y}+\|\bar{v}\|_{L^\infty}+\|\bar{v}\|_{W^{1,4}})^3\le C(1+\sqrt{Y_*}+C(1+Y_*))^3= M_*.\]
	This gives the second bound in \eqref{est:bound}.
\end{proof}

\begin{remark}\label{rem:regularity_n}
As a consequence, we additionally obtain from \eqref{est:y} that
\[\e\int_0^{T_\e}(\|\nabla w\|_{L^2}^2+\|\Delta w\|_{L^2}^2)\,dt\le C(Y_*,T_\e),\]
which implies, by the elliptic estimate $\|w\|_{H^2}\le C(\|\Delta w\|_{L^2}+\|w\|_{L^2})$ for $w$ with $\pa_\nu w=0$, $\|\pa_tn\|_{L^2H^2}=\|w\|_{L^2H^2}\le C(\e,Y_*)$. Moreover, it follows from the definition of $\tau=\tau_{\bar{v}}(n)$ and Lemma \ref{lem:tau} (3) that
\[\Delta n =\tau-|\nabla n|^2n-n\times (\bar{v}\cdot\nabla)n,\quad \mbox{with}\quad \tau =\frac{\e w-n\times w}{1+\e^2}.\]
Thus, since $H^2$ is an algebra and $\|\bar{v}\|_{H^2}\le C\|\bar{n}\|_{H^2}\|\bar{n}\|_{H^3}\le C(M_*)$ by Lemma \ref{lem:velocity_map} (2), we get
\begin{align*}
	\|\Delta n\|_{H^2}&\le \|\tau\|_{H^2} +C\|n\|_{H^3}^2\|n\|_{H^2}+\|n\|_{H^2}\|\bar{v}\|_{H^2}\|n\|_{H^3}\\
	&\le C\|w\|_{H^2}+C\|n\|_{H^2}\|w\|_{H^2}+C\|n\|_{H^3}^2\|n\|_{H^2}+\|n\|_{H^2}\|\bar{v}\|_{H^2}\|n\|_{H^3}\\
	&\le C(M_*)(1+\|w\|_{H^2}),
\end{align*}
from which we conclude that $\|n\|_{L^2H^4}\le C(\e,Y_*)$.	
\end{remark}

Therefore, Proposition \ref{prop:invariance} and Remark \ref{rem:regularity_n} imply that there exist constants $Y_*$, $M_*$, and $T_\e$ such that we can construct inductively the sequence of approximated solutions $n^{k}\in\mathcal{C}_{T_\e}$ via \eqref{eq:iteration} that satisfies
\[\sup_{0\le t\le T_\e}(\|\Delta n^{k}(t)\|_{L^2}^2+\|\pa_tn^{k}(t)\|_{H^1}^2)\le Y_*,\quad \sup_{0\le t\le T_\e}\|n^{k}(t)\|_{H^3}\le M_*,\]
and
\[\|n^{k}\|_{L^2H^4}\le C(\e,Y_*),\quad \|\pa_t n^{k}\|_{L^2H^2}\le C(\e,Y_*).\]
Here, for the first step $k=0$, $\bar{n}=n^{0}=n_0$ is time-independent. The proof of Proposition \ref{prop:invariance} uses $\bar{n}$ only through $\bar{v}=K[n_0]$, which belongs to $\mathcal{V}_{T_\e}$ with $\pa_t\bar{v}=0$ by Lemma \ref{lem:velocity_map}, and the bound $\|n_0\|_{H^3}\le M_*$ holds by applying Lemma \ref{lem:norm_equiv} to $n^{1}$ at $t=0$, exactly as in the last step of the proof of Proposition \ref{prop:invariance}.

\subsection{Contraction in low-order norm}
Next, we show that the sequence $\{n^{k}\}_{k=1}^\infty$ is a Cauchy sequence in a low-order function space, after we choose further small time $T_\e$. 

\begin{lemma}\label{lem:cauchy}
	After taking $T_\e$ further small if needed, the sequence $\{n^{k}\}_{k=1}^\infty$ is a Cauchy sequence in $C([0,T_\e];L^2)$.
\end{lemma}
\begin{proof}
	Define $z^{k+1}:=n^{k+1}-n^{k}$ and $q^{k}:=v^{k}-v^{k-1}$. Subtracting the equations for $n^{k+1}$ and $n^{k}$, $z^{k+1}$ satisfies
	\begin{align}
		\begin{aligned}\label{eq:z}
		\pa_t z^{k+1} &= \e\Delta z^{k+1}+\e\left(|\nabla n^{k+1}|^2n^{k+1}-|\nabla n^{k}|^2n^{k}\right)\\
		&\quad+\e(z^{k+1}\times(v^{k}\cdot\nabla n^{k+1})+n^{k}\times(v^{k}\cdot\nabla z^{k+1})+n^{k}\times(q^{k}\cdot\nabla n^{k}))\\
		&\quad + z^{k+1}\times\Delta n^{k+1}+n^{k}\times\Delta z^{k+1}-(v^{k}\cdot\nabla )z^{k+1}-(q^{k}\cdot\nabla)n^{k}.
		\end{aligned}
	\end{align}
	By testing \eqref{eq:z} with $z^{k+1}$, we get
	\begin{align*}
		\frac{1}{2}\frac{d}{dt}\|z^{k+1}\|^2 &= \e\<z^{k+1},\Delta z^{k+1}\>+\e \<z^{k+1},|\nabla n^{k+1}|^2n^{k+1}-|\nabla n^{k}|^2n^{k}\>\\
		&\quad +\e \<z^{k+1},z^{k+1}\times(v^{k}\cdot\nabla n^{k+1})+n^{k}\times(v^{k}\cdot\nabla z^{k+1})+n^{k}\times(q^{k}\cdot\nabla) n^{k}\>\\
		&\quad + \<z^{k+1},z^{k+1}\times\Delta n^{k+1}+n^{k}\times\Delta z^{k+1}-(v^{k}\cdot\nabla )z^{k+1}-(q^{k}\cdot\nabla)n^{k}\>\\
		&=:-\e\|\nabla z^{k+1}\|_{L^2}^2+\sum_{i=1}^3 K_i,
	\end{align*}
	where $\<\cdot,\cdot\>$ denotes the $L^2$ inner product.\\
	
	\noindent $\bullet$ (Estimate of $K_1$): We first split $K_1$ as
	\begin{align*}
		K_1 = \e\<z^{k+1},(|\nabla n^{k+1}|^2-|\nabla n^k|^2)n^{k+1}\>+\e\<z^{k+1},|\nabla n^k|^2z^{k+1}\>.
	\end{align*}
	Since $|n^{k+1}|=|n^k|=1$, we have $|z^{k+1}|^2=2-2n^{k+1}\cdot n^k$. This, together with $|z^{k+1}|\le 2$, yields
	\[z^{k+1}\cdot n^{k+1} = 1-n^{k+1}\cdot n^k=\frac{1}{2}|z^{k+1}|^2\le |z^{k+1}|.\]
	Therefore, the first term of $K_1$ can be bounded as
	\begin{align*}
		\<z^{k+1},(|\nabla n^{k+1}|^2-|\nabla n^k|^2)n^{k+1}\>&=\<z^{k+1},(\nabla(n^{k+1}+n^k):\nabla z^{k+1})n^{k+1}\>\\
		& \le (\|\nabla n^{k+1}\|_{L^\infty}+\|\nabla n^k\|_{L^\infty})\|\nabla z^{k+1}\|_{L^2}\|z^{k+1}\|_{L^2}.
	\end{align*}
	This yields
	\[K_1\le \e(\|\nabla n^{k+1}\|_{L^\infty}+\|\nabla n^k\|_{L^\infty})\|\nabla z^{k+1}\|_{L^2}\|z^{k+1}\|_{L^2}+\e\|\nabla n^k\|_{L^\infty}^2\|z^{k+1}\|_{L^2}^2.\]
	
	\noindent $\bullet$ (Estimate of $K_2$): Since $z^{k+1}\cdot(z^{k+1}\times(v^k\cdot\nabla n^{k+1}))=0$,  we estimate $K_2$ as
	\begin{align*}
		K_2 &= 0+ \e \<z^{k+1},n^k\times(v^k\cdot\nabla z^{k+1})\>+\e\<z^{k+1},n^k\times (q^k\cdot\nabla n^k)\>\\
		&\le \|v^k\|_{L^\infty}\|\nabla z^{k+1}\|_{L^2}\|z^{k+1}\|_{L^2}+\|q^k\|_{L^2}\|\nabla n^{k}\|_{L^\infty}\|z^{k+1}\|_{L^2}.
	\end{align*}
	\noindent $\bullet$ (Estimate of $K_3$): Similarly, we separate $K_3$ as
	\begin{align*}
		K_3 = 0+\<z^{k+1},n^k\times \Delta z^{k+1}\>-\<z^{k+1},v^k\cdot\nabla z^{k+1}\>-\<z^{k+1},q^k\cdot\nabla n^k\>.
	\end{align*}
	Since
	\[\<z^{k+1},n^k\times\Delta z^{k+1}\>=\<\Delta z^{k+1},z^{k+1}\times n^k\>=-\<\pa_j z^{k+1},z^{k+1}\times\pa_j n^k\>,\]
	and
	\[\<z^{k+1},v^k\cdot\nabla z^{k+1}\>=0,\]
	we get
	\[K_3\le \|\nabla n^k\|_{L^\infty}\|z^{k+1}\|\|\nabla z^{k+1}\|+\|q^k\|_{L^2}\|\nabla n^k\|_{L^\infty}\|z^{k+1}\|_{L^2}.\]
	Combining the estimates for $K_i$ for $i=1,2,3$, and using the uniform bounds
	\[\|\nabla n^{k+1}\|_{L^\infty},\|\nabla n^k\|_{L^\infty}, \|v^k\|_{L^\infty}\le C(M_*),\] 
	we get
	\begin{align*}
		\frac{1}{2}\frac{d}{dt}&\|z^{k+1}\|_{L^2}^2+\e\|\nabla z^{k+1}\|_{L^2}^2\\
		&\le C(M_*)\|\nabla z^{k+1}\|_{L^2}\|z^{k+1}\|_{L^2}+C(M_*)\|z^{k+1}\|_{L^2}^2+C(M_*)\|q^k\|_{L^2}\|z^{k+1}\|_{L^2}\\
		&\le K_\e\|z^{k+1}\|_{L^2}^2+\|q^k\|_{L^2}^2+\frac{\e}{2}\|\nabla z^{k+1}\|^2_{L^2},
	\end{align*}
	for some positive constant $K_\e$ that may depend on $M_*$ and $\e$. Now, it follows from Lemma \ref{lem:velocity_map} (3) that 
	\[\|q^k\|_{L^2}\le C(M_*)\|z^k\|_{L^2}, \]
	which implies that 
	\[\frac{1}{2}\frac{d}{dt}\|z^{k+1}\|_{L^2}^2+\frac{\e}{2}\|\nabla z^{k+1}\|_{L^2}^2\le K_\e\|z^{k+1}\|_{L^2}^2+C(M_*)\|z^k\|_{L^2}^2,\quad z^{k+1}(0)=0.\]
	Then, the Gr\"onwall inequality implies
	\[\sup_{0\le t\le T_\e}\|z^{k+1}\|_{L^2}^2\le C(M_*)T_\e e^{K_\e T_\e}\sup_{0\le t\le T_\e}\|z^k\|_{L^2}^2.\]
	Thus, if we further take $T_\e$ sufficiently small so that
	\[C(M_*)T_\e e^{K_\e T_\e}\le \frac{1}{4},\]
	we get
	\[\sup_{0\le t\le T_\e}\|z^{k+1}\|_{L^2}\le \frac{1}{2}\sup_{0\le t\le T_\e}\|z^k\|_{L^2},\]
	and therefore, the sequence $\{n^k\}_{k=1}^\infty$ is Cauchy in $C([0,T_\e];L^2)$. 
\end{proof}

\subsection{Proof of Proposition \ref{prop:regularize_exist}}

As a consequence of Lemma \ref{lem:cauchy}, there exists a strong limit $n=\lim_{k\to\infty}n^{k}$ that belongs to $C([0,T_\e];L^2)$. Furthermore, since the sequence $\{n^k\}_{k=1}^\infty$ is uniformly bounded in $L^\infty H^3\cap L^2H^4$ and the sequence $\{\pa_tn^k\}_{k=1}^\infty$ is uniformly bounded in $L^\infty H^1\cap L^2H^2$, the limit $n$ also belongs to $L^\infty H^3\cap L^2H^4$ with $\pa_tn \in L^\infty H^1\cap L^2H^2$. Then, using Lions--Magenes lemma, we conclude that $n\in C([0,T_\e];H^3)$. In particular, we have $n(0)=n_0$, since $n^{k}(0)=n_0$ for all $k$. 
Now, we pass to the limit $k\to\infty$ to show that $n$ indeed solves the regularized equation \eqref{eq:ISF_regularize}. We first observe that 
\[n^k\to n \quad \mbox{in}\quad C([0,T_\e]; H^{3-\delta})\]
for any $0<\delta<3$ by interpolating convergence in $L^2$ and uniform boundedness in $H^3$, and after taking the subsequence, we may assume that $n^k\to n$ pointwisely. Then, by Lemma \ref{lem:velocity_map} (3), we have
\[v^{k}\to v:=K[n]\quad\mbox{in}\quad C([0,T_\e];H^1).\]
Furthermore, since
\begin{align*}
	n^k&\cdot(\pa_1 n^k\times \pa_2 n^k)-n\cdot(\pa_1 n\times\pa_2n)\\
	&=(n^k-n)\cdot (\pa_1 n^k\times \pa_2 n^k)+ n\cdot (\pa_1 (n^k-n)\times \pa_2 n^k)+n\cdot (\pa_1 n\times \pa_2 (n^k-n)),
\end{align*}
and $H^{2-\delta}$ is an algebra, we have
\[\|n^k\cdot(\pa_1 n^k\times \pa_2 n^k)-n\cdot(\pa_1 n\times\pa_2n)\|_{H^{2-\delta}}\le C(M_*)\|n^k-n\|_{H^{3-\delta}}\to0,\]
from which the convergence of the velocity field follows:
\[\|v^{k}-v\|_{H^{3-\delta}}\le C\|\omega^k-\omega\|_{H^{2-\delta}}\to 0,\quad\mbox{where}\quad \omega^k:=n^k\cdot(\pa_1 n^k\times \pa_2 n^k),\quad \omega:=n\cdot(\pa_1 n\times\pa_2n).\]
Using the above convergences, we now take the limit in the equation \eqref{eq:ISF_regularize}:
\[\pa_t n^{k+1}=\e\tau^k+n^{k+1}\times \tau^k,\quad \tau^k = \Delta n^{k+1}+|\nabla n^{k+1}|^2n^{k+1}+n^{k+1}\times (v^k\cdot\nabla)n^{k+1}.\]
From the above strong convergences, it is straightforward to obtain
\begin{align*}
	&\Delta n^{k+1}\to\Delta n\quad\mbox{in}\quad H^{1-\delta}\\
	&|\nabla n^{k+1}|^2n^{k+1}\to|\nabla n|^2n\quad\mbox{in}\quad H^{2-\delta}\\
	&n^{k+1}\times(v^k\cdot\nabla)n^{k+1}\to n\times (v\cdot\nabla)n \quad\mbox{in}\quad H^{2-\delta}.
\end{align*}
Hence, we get $\tau^k\to\tau$ in $C([0,T_\e];H^{1-\delta})$, which implies
\[\|n^{k+1}\times \tau^k-n\times \tau\|_{L^2}\le \|n^{k+1}-n\|_{L^\infty}\|\tau^k\|_{L^2}+\|n\|_{L^\infty}\|\tau^k-\tau\|_{L^2}\to0.\]
The unit length condition $|n(t,x)|=1$ for all $(t,x)\in [0,T_\e]\times \Omega$ holds since $n$ is a pointwise limit of $n^k$, and the Neumann boundary condition $\pa_\nu n=0$ on $\pa\Omega$ also holds since $n^k\to n$ in $C([0,T_\e];H^{3-\delta})$ and the normal trace is a bounded linear map from $H^{3-\delta}(\Omega)$ to $L^2(\pa\Omega)$. Also, since $\pa_tn^{k+1}\to\pa_tn$ in the sense of distributions as $n^{k+1}\to n$ in $C([0,T_\e];L^2)$, we take $k\to\infty$ to observe that $n$ solves the regularized equation \eqref{eq:ISF_regularize} in the sense of distributions on $(0,T_\e)\times\Omega$, and hence almost everywhere since both sides of the equation belong to $L^\infty(0,T_\e;L^2)$. 

Furthermore, since the right-hand side of the equation for $w=\pa_tn$ in \eqref{eq:w-0} belongs to $L^2L^2$, we conclude that $\pa_tw=\pa_{tt}n\in L^2L^2$ and again using Lions--Magenes lemma, we have $\pa_tn\in C([0,T_\e];H^1)$. This shows that the limit $n$ has the desired regularity stated in Proposition \ref{prop:regularize_exist}.

Finally, to show the uniqueness of the solution \eqref{eq:ISF_regularize}, let $n,\tilde{n}\in\mathcal{C}_{T_\e}$ be two solutions of \eqref{eq:ISF_regularize} with the same initial data $n_0$, and let $M=\max\{\|n\|_{L^\infty H^3},\|\tilde{n}\|_{L^\infty H^3}\}$. Then, it follows from the similar computation as we did in Lemma \ref{lem:cauchy} that the difference $z:=n-\tilde{n}$ satisfies
\[\frac{1}{2}\frac{d}{dt}\|z\|_{L^2}^2+\frac{\e}{2}\|\nabla z\|_{L^2}^2\le C(M,\e)\|z\|_{L^2}^2,\quad z(0)=0,\]
from which we deduce the uniqueness by applying the Gr\"onwall inequality. This completes the proof of Proposition \ref{prop:regularize_exist}.

\section{Proof of Theorem \ref{thm:director}}\label{sec:4}
\setcounter{equation}{0}

In this section, we obtain the uniform-in-$\e$ estimates for the solution $n_\e$ to \eqref{eq:ISF_regularize} constructed in the previous section. After obtaining uniform bound, we extract the subsequence as $\e\to0$ whose limit satisfies the equation \eqref{eq:ISF_director_h2}, completing the proof of Theorem \ref{thm:director}. We also show that the solution to \eqref{eq:ISF_director_h2} is unique.

\subsection{Uniform-in-$\e$ estimate and existence of the solution to \eqref{eq:ISF_director_h2}}
Up to now, we have constructed the solution $n_\e\in \mathcal{C}_{T_\e}$ to \eqref{eq:ISF_regularize} for each fixed $\e$, whose existence time may depend on $\e$. Below, we show that, these solutions $n_\e$ has a common existence time $T_*$, which is independent of $\e$.

\begin{proposition}
	There exists $T_*$ and $M>0$, depending only on the initial data $n_0$ such that for all $\e\in(0,1)$, the solution $n_\e\in\mathcal{C}_{T_\e}$ of \eqref{eq:ISF_regularize} constructed in Proposition \ref{prop:regularize_exist} can be extended to $[0,T_*]$ and it satisfies
	\begin{equation}\label{uniform_bound_e}
		\sup_{0\le t\le T_*}\left(\|n_\e\|_{H^3}^2+\|\pa_tn_\e\|_{H^1}^2\right)\le M.
	\end{equation}
\end{proposition}

\begin{proof}
	For simplicity, we drop the dependency of $\e$ in the proof, that is $n=n_\e$. Define $y$ and $\Gamma$ as 
	\[y(t):=\|\Delta n(t)\|_{L^2}^2+\|\pa_tn(t)\|_{H^1}^2,\]
	and
	\[\Gamma(t):=1+\|n(t)\|_{H^2}+\|\pa_tn(t)\|_{H^1}+\|v(t)\|_{L^\infty}+\|v(t)\|_{W^{1,4}}.\]
	In particular, the velocity $v$ in $\Gamma$ is now generated from $n$, that is $v=K[n]$. Then, following exactly the same computation as we did in Proposition \ref{prop:invariance} to derive \eqref{est:dydt}, we get
	\begin{equation}\label{dydt}
		\frac{dy}{dt}\le C(1+y+\|v\|_{L^\infty}^2+\|v\|_{W^{1,4}}^2)^4+C(\|\nabla v\|_{L^\infty}+\|\pa_t v\|_{H^1}^2)(1+y+\|v\|_{L^\infty}^2+\|v\|_{W^{1,4}}^2).
	\end{equation}
	Furthermore, we use Lemma \ref{lem:velocity_map} (2) to get
	\[\|v\|_{L^\infty}+\|v\|_{W^{1,4}}\le C\|n\|_{H^2}^2\le C(1+y),\]
	from which we bound $\Gamma$ in terms of $y$ as
	\[\Gamma\le C(1+\sqrt{y}+y)\le C(1+y).\]
	Furthermore, since $n$ satisfies \eqref{eq:ISF_regularize} with $v=K[n]\in\mathcal{V}_T$, Lemma \ref{lem:norm_equiv} is valid with the velocity $v$. This implies $\|n\|_{H^3}\le C\Gamma^3$ and therefore,  
	\[\|\pa_tv\|_{H^1}\le C\|n\|_{H^3}\|\pa_tn\|_{H^1}\le C(1+y)^3\sqrt{y},\]
	and 
	\[\|\nabla v\|_{L^\infty}\le C\|n\|_{H^2}\|n\|_{H^3}\le C(1+y)^4.\]
	Substituting the above estimates into \eqref{dydt}, we observe that $y$ satisfies
	\[\frac{dy}{dt}\le C(1+y)^8+C((1+y)^4+(1+y)^7)(1+y)^2\le C_0(1+y)^9,\]
	which implies
	\[1+y(t)\le (1+y_0)\left(1-8C_0(1+y_0)^8t\right)^{-1/8}\le 2(1+y_0),\quad\mbox{for}\quad t\le \frac{1-2^{-8}}{8C_0(1+y_0)^8}.\]
	Recall that the initial value $y_0$ satisfies for $\e<1$ that
	\begin{align*}
		y_0 &= \|\Delta n_0\|_{L^2}^2 + \|w_0\|_{H^1}^2=\|\Delta n_0\|_{L^2}^2 + \|\e\tau_{v_0}(n_0) + n_0\times \tau_{v_0}(n_0)\|_{H^1}^2\\
		&\le \|\Delta n_0\|_{L^2}^2 + C(1+\|n_0\|^2_{H^3})\|\tau_{v_0}(n_0)\|^2_{H^1}=:y^*_0.
	\end{align*}
	Accordingly, we define $\e$-independent time $T_*$ as
	\[T_*:=\frac{1-2^{-8}}{8C_0(1+y^*_0)^8}.\] 
	Since $y_0\le y_0^*$, the bound $1+y(t)\le 2(1+y_0)$ holds for all $t\le T_*$. Now, we show that the solution $n=n_\e\in \mathcal{C}_{T_\e}$ can be extended to $T_*$. Suppose the existence time $T_{\max}$ of $n$ is less than $T_*$. Then, on the interval $[0,T_{\max}]$, we have
	\[y(t)\le 1+2y_0\le Y_*,\quad\mbox{and}\quad \|n\|_{H^3}\le C(1+y)^3\le C(3+2y^*_0)^3=:M_1.\]
	
	On the other hand, it follows from Lemma \ref{lem:velocity_map} (2) and the Sobolev embedding that
	\[\|\tau_{v_0}(n_0)\|_{H^1}\le C(1+\|n_0\|_{H^3})^4.\]
	Therefore, $y_0^*$, and hence the existence time $T_\e$ of Proposition \ref{prop:regularize_exist}, can be bounded in terms of $\|n_0\|_{H^3}$. This implies that the existence time provided by Proposition \ref{prop:regularize_exist} for the initial data $n(t_1)$ for $t_1\in[0,T_{\max})$ with $\|n(t_1)\|_{H^3}\le M_1$ can be taken to depend only on $M_1$ and $\e$. We denote it by $T_\e(M_1)$. Then, for any $t_1\in [0,T_{\max})$, we apply Proposition \ref{prop:regularize_exist} with the initial data $n(t_1)$, and the solution exists in the time interval $[t_1,t_1+T_\e(M_1)]$. Since $t_1\in[0,T_{\max})$ is arbitrary, one may choose $t_1$ so that $t_1>T_{\max}-T_\e(M_1)$. Then, the solution can be extended until $t_1+T_\e(M_1)>T_{\max}$, which contradicts the maximality of the existence time. Therefore, for each fixed $\e$, the solution $n=n_\e$ can be extended up to time $T_*$, which is $\e$-independent, and \eqref{uniform_bound_e} holds with $M:=M_1^2+1+2y_0^*$.
\end{proof}

Since the solution $n_\e$ is guaranteed to exist on the common time interval $[0,T_*]$, we are now able to take the limit $\e\to0$. Thanks to \eqref{uniform_bound_e}, the sequence $\{n_\e\}_{\e>0}$ is bounded in $L^\infty(0,T_*;H^3)$ and $\{\pa_tn_\e\}_{\e>0}$ is bounded in $L^\infty H^1$. Therefore, we apply Aubin--Lions lemma to conclude that there exists a sequence $\e_k\to0$ such that $n_{\e_k}\to n$ in $C([0,T_*];H^{3-\delta})$ with $\pa_tn_{\e_k}\xrightharpoonup{*}\pa_t n$ in $L^\infty H^1$, and $n\in L^\infty H^3$. Furthermore, by extracting further subsequence, we also obtain that $n_{\e_k}\to n$ pointwisely. Hence, by using the same argument as in the previous section, we have the following convergences:
	\begin{align*}
		v_{\e_k}\to v:=K[n] \quad\mbox{in}\quad  C([0,T_*];H^{3-\delta}),\\
		n_{\e_k}\times \tau_{\e_k} \to n\times \tau \quad\mbox{in}\quad C([0,T_*];L^2).
	\end{align*}
Furthermore, since each terms in $\tau_\e=\Delta n_\e+|\nabla n_\e|^2n_\e+n_\e\times(v_\e\cdot\nabla)n_\e$ are bounded in $L^\infty L^2$, we get $\e_k\tau_{\e_k}\to0$ in $L^\infty L^2$ as $\e_k\to0$. Moreover, $\pa_tn_{\e_k}\to\pa_tn$ in the sense of distributions, so that the limit $n$ satisfies \eqref{eq:ISF_director_h2} in the sense of distributions on $(0,T_*)\times\Omega$, and hence almost everywhere, since both sides of the equation belong to $L^\infty(0,T_*;L^2)$. The unit length condition $|n(t,x)|=1$ and the Neumann boundary condition $\pa_\nu n=0$ on $\pa\Omega$ also hold as in the proof of Proposition \ref{prop:regularize_exist}, and $n(0)=n_0$, since $n_{\e_k}(0)=n_0$ for all $k$ and $n_{\e_k}\to n$ in $C([0,T_*];H^{3-\delta})$. Finally, the bound \eqref{uniform_bound_e} is inherited by $n$ by the lower semicontinuity of the norms under the weak-$*$ convergence. This shows that the limit $n$ is a strong solution to \eqref{eq:ISF_director_h2} in the sense of Definition \ref{def:director}, satisfying the bound stated in Theorem \ref{thm:director}. 

\subsection{Uniqueness of the solution to \eqref{eq:ISF_director_h2}}
We now show the uniqueness of the solution to \eqref{eq:ISF_director_h2}.

\begin{lemma}\label{lem:uniqueness}
	Let $n,\tilde{n}$ be two strong solutions to \eqref{eq:ISF_director_h2} on $[0,T_*]$ in the sense of Definition \ref{def:director} with the common initial data $n_0$. Then $n=\tilde{n}$.
\end{lemma}
\begin{proof}
	Let $M$ be a constant satisfying $\|n\|_{L^\infty H^3},\|\tilde{n}\|_{L^\infty H^3}<M$, and define
	\[z:=n-\tilde{n},\quad q:=v-\tilde{v}:=K[n]-K[\tilde{n}].\]
	Then, subtracting the equations for $n$ and $\tilde{n}$, we get
	\begin{equation}\label{eq:z_lim}
		\pa_t z = z\times \Delta n +\tilde{n}\times \Delta z-(v\cdot\nabla)z -(q\cdot\nabla)\tilde{n}.
	\end{equation}
	We multiply \eqref{eq:z_lim} with $z$ and use the incompressibility of $v$ to get
	\begin{align*}
		\frac{1}{2}\frac{d}{dt}\|z\|^2&=0+\int_{\Omega}z\cdot(\tilde{n}\times\Delta z)\,dx - 0 - \int_{\Omega} z\cdot (q\cdot\nabla) \tilde{n}\,dx\\
		&=-\int_{\Omega} \pa_j z\cdot(z\times \pa_j\tilde{n})\,dx-\int_{\Omega}z\cdot(q\cdot\nabla)\tilde{n}\,dx\\
		&\le C\|\nabla\tilde{n}\|_{L^\infty}\|z\|(\|\nabla z\|+\|q\|)\le C(M)\|z\|_{H^1}^2.
	\end{align*}
	Next, we estimate the first order derivative of $z$ by using \eqref{eq:z_lim} as
	\begin{align*}
		\frac{1}{2}\frac{d}{dt}\|\nabla z\|_{L^2}^2&=-\int_\Omega \Delta z\cdot\pa_t z\,dx\\
		&=-\int_\Omega \Delta z\cdot( z\times\Delta n)\,dx -\int_{\Omega}\Delta z\cdot(\tilde{n}\times \Delta z)\,dx \\
		&\quad +\int_\Omega \Delta z\cdot (v\cdot\nabla) z\,dx +\int_{\Omega}\Delta z\cdot (q\cdot\nabla)\tilde{n}\,dx\\
		&=\int_\Omega \pa_j z\cdot(z\times \pa_j\Delta n)\,dx-0-\int_\Omega \pa_jv_i\pa_iz\cdot \pa_j z\,dx-\int_{\Omega}\pa_j z\cdot \pa_j((q\cdot\nabla)\tilde{n})\,dx\\
		&\le \|\nabla z\|_{L^2}\|z\|_{L^\infty}\|\nabla\Delta n\|_{L^2}+\|\nabla v\|_{L^\infty}\|\nabla z\|_{L^2}^2\\
		&\quad +\|\nabla z\|_{L^2}(\|\nabla q\|\|\nabla \tilde{n}\|_{L^\infty} +\|q\|_{L^4}\|\nabla^2\tilde{n}\|_{L^4})\\
		&\le C(M)\|z\|_{H^1}\|z\|_{L^\infty}+C(M)\|z\|_{H^1}^2.
	\end{align*}
	Now, we use the Brezis--Gallou\"et inequality \cite{BG80}
	\[\|f\|_{L^\infty}\le C\|f\|_{H^1}\left(1+\log \left(1+\frac{\|f\|_{H^2}}{\|f\|_{H^1}}\right)\right)^{1/2},\]
	elementary inequalities 
	\[(1+\log (1+a))^{1/2}\le 1+\log(1+a)\le 1+\log 2\left(1+a^2\right),\quad \mbox{for}\quad a>0,\]
	and $\|z\|_{H^2}\le 2M$ to conclude that
	\begin{align*}
		\frac{d}{dt}\|z\|_{H^1}^2 &\le C(M)\|z\|_{H^1}^2+C(M)\|z\|_{H^1}^2\left(1+\log\left(1+\frac{2M}{\|z\|_{H^1}}\right)\right)^{1/2}\\
		&\le C(M)\|z\|_{H^1}^2+C(M)\|z\|_{H^1}^2\left(1+\log 2+ \log \left(1+\frac{4M^2}{\|z\|_{H^1}^2}\right)\right)\\
		&\le C(M)\|z\|_{H^1}^2\left(1+\log\left(1+\frac{4M^2}{\|z\|_{H^1}^2}\right)\right)=:f(\|z\|_{H^1}^2).
	\end{align*}
	Since the function $f(r)=C(M)r(1+\log(1+\frac{4M^2}{r}))$ is positive increasing function satisfying
	\[\int_0^1 \frac{1}{f(s)}\,ds=+\infty,\]
	it satisfies the condition for the Osgood lemma (see, e.g., \cite[Lemma 3.4]{BCD11}). Hence, we conclude that $\|z\|_{H^1}=0$, which proves the uniqueness of the solution.
\end{proof}

Combining the existence of the solution $n$ to \eqref{eq:ISF_director_h2} and Lemma \ref{lem:uniqueness}, we complete the proof of Theorem \ref{thm:director}.

\begin{remark}
	As we did in Section \ref{sec:3}, we scaled the time variable as $t\mapsto\frac{\hbar}{2}t$ to transform \eqref{eq:ISF_director} into \eqref{eq:ISF_director_h2}. Since we proved that the $\hbar$-independent system \eqref{eq:ISF_director_h2} has a unique local solution up to time $T_*$ for an $O(1)$ initial data $n_0$, the \eqref{eq:ISF_director} has a unique local solution up to time $\frac{2T_*}{\hbar}$, which diverges to $+\infty$ in the semi-classical limit regime $\hbar\to0$.
\end{remark}

\section{Proof of Theorem \ref{thm:spinor}}\label{sec:5}
\setcounter{equation}{0}
In this section, we prove Theorem \ref{thm:spinor} by constructing the solution $\psi$ to \eqref{eq:ISF} from the solution $n$ to \eqref{eq:ISF_director}, whose existence and uniqueness is guaranteed by Theorem \ref{thm:director}. 

\subsection{Construction of $\psi$}
From the given initial data $\psi_0\in H^3$, we construct $n_0=\psi_0^\dg\sigma\psi_0\in H^3$. Then, by Theorem \ref{thm:director}, there exists a positive time $T_*$ and the unique solution $n\in L^\infty(0,T_*;H^3)$. In this part, we construct $\psi$ from $n$, which is a possible candidate of the solution to \eqref{eq:ISF}. To this end, for any two unit vectors $n,n_0\in\bbs^2$ with $n\cdot n_0>-1$, we define a matrix $U(n,n_0)$ as
\[U(n,n_0):=\frac{(1+n\cdot n_0)I+\i (n\times n_0)\cdot\sigma}{\sqrt{2(1+n\cdot n_0)}}.\]
The motivation of introducing the unitary matrix $U(n,n_0)$ is to evolve the initial spinor $\psi_0$ to $\tilde{\psi} = U(n,n_0)\psi_0$ that satisfies $\tilde{\psi}^\dg \sigma \tilde{\psi}=n$, as we will verify soon. 

\begin{lemma}\label{lem:U}
	For any given $n,n_0\in \bbs^2$, with $n\cdot n_0>-1$, the matrix $U(n,n_0)$ is unitary, and it satisfies
	\[U(n,n_0)(n_0\cdot\sigma)U(n,n_0)^\dg = n\cdot\sigma.\]
\end{lemma}
\begin{proof}
Since $\sigma_j$ are Hermitian, 
\begin{align*}
	U(n,n_0)^\dg U(n,n_0)&=\frac{1}{2(1+n\cdot n_0)}\left((1+n\cdot n_0)^2I + ((n\times n_0)\cdot \sigma)((n\times n_0)\cdot\sigma)\right)\\
	&=\frac{1}{2(1+n\cdot n_0)}\left((1+n\cdot n_0)^2I + |n\times n_0|^2I\right)=I,
\end{align*}
where we used \eqref{Pauli_identity_2} and
\[(1+n\cdot n_0)^2+|n\times n_0|^2=1+2n\cdot n_0+|n|^2|n_0|^2=2(1+n\cdot n_0).\]
Furthermore, direct computation gives 
\begin{align*}
	&U(n,n_0) (n_0\cdot\sigma) U(n,n_0)^\dg\\
	&=\frac{1}{2(1+n\cdot n_0)}\left((1+n\cdot n_0)I+\i(n\times n_0)\cdot \sigma\right)(n_0\cdot \sigma)\left((1+n\cdot n_0)I-\i(n\times n_0)\cdot \sigma\right)\\
	&=\frac{1}{2(1+n\cdot n_0)}\left((1+n\cdot n_0)(n_0\cdot \sigma)-((n\cdot n_0)n_0-n)\cdot\sigma \right)\left((1+n\cdot n_0)I-\i(n\times n_0)\cdot \sigma\right)\\
	&=\frac{1}{2(1+n\cdot n_0)}\left(n_0\cdot \sigma+n\cdot\sigma \right)\left((1+n\cdot n_0)I-\i(n\times n_0)\cdot \sigma\right)\\
	&=\frac{1}{2(1+n\cdot n_0)}\left((1+n\cdot n_0)(n_0\cdot \sigma)+(1+n\cdot n_0)(n\cdot\sigma) \right)\\
	&\quad +\frac{1}{2(1+n\cdot n_0)}\left((n-(n_0\cdot n)n_0)\cdot\sigma+((n\cdot n_0)n-n_0)\cdot\sigma\right)\\
	&=n\cdot\sigma.
\end{align*}
\end{proof}

Since the solution $n$ to \eqref{eq:ISF_director} belongs to $C([0,T_*];H^{3-\delta})\subset C([0,T_*];C(\bar{\Omega}))$, we have $\|n(t)-n_0\|_{L^\infty}\to0$ as $t\to0$, and since $|n_0(x)|^2=1$ for all $x\in\Omega$, there exists $T_{**}\in(0,T_*]$ such that $n_0(x)\cdot n(t,x)>-1$ for all $x\in \Omega$ and $t\in[0,T_{**}]$. Therefore, once we define $\tilde{\psi}(t,x):=U(n(t,x),n_0(x))\psi_0(x)$ for $t\in[0,T_{**}]$, we use Lemma \ref{lem:U} to observe that
\begin{align*}
	\tilde{\psi}^\dg\sigma_k\tilde{\psi}&=(U\psi_0)^\dg \sigma_k (U\psi_0)=\tr(\sigma_k(U\psi_0)(U\psi_0)^\dg)=\tr(\sigma_k U\psi_0\psi_0^\dg U^\dg)\\
	&=\frac{1}{2}\tr(\sigma_k U(I+n_0\cdot\sigma) U^\dg)=\frac{1}{2}\tr(\sigma_k +\sigma_k(n\cdot\sigma))=n_k.
\end{align*}
Therefore, the spinor $\tilde{\psi}$ exactly generates the director field $n$. Unfortunately, this constructed $\tilde{\psi}$ needs not satisfy the incompressibility condition \[\dv\,(\hbar\Im(\tilde{\psi}^\dg\nabla\tilde{\psi}))=\hbar\Im(\tilde{\psi}^\dg\Delta\tilde{\psi})=0\]
in general. To attain incompressibility, we consider the phase-corrected spinor $\psi = e^{\i\theta}\tilde{\psi}$, where $\theta=\theta(t,x)$ is the solution to the following Poisson equation with the Neumann boundary condition, normalized by $\int_\Omega\theta\,dx=0$:
\begin{equation}\label{eq:poisson}
	\hbar\Delta\theta = -\dv\,(\hbar\Im(\tilde{\psi}^\dg\nabla\tilde{\psi})),\quad \pa_\nu\theta=0.
\end{equation}
Note that since
\[\int_\Omega \dv (\Im(\tilde{\psi}^\dg\nabla\tilde{\psi}))\,dx=\int_{\pa\Omega}\Im((\tilde{\psi}^\dg\pa_\nu\tilde{\psi}))\,ds=0,\]
the Neumann problem \eqref{eq:poisson} is solvable. Then, we observe that
\[\psi^\dg\pa_j\psi = e^{-\i\theta}\tilde{\psi}^\dg (e^{\i\theta}\i\pa_j\theta\tilde{\psi}+e^{\i\theta}\pa_j\tilde{\psi})=\i\pa_j\theta +\tilde{\psi}^\dg\pa_j\tilde{\psi}.\]
Hence, multiplying $\hbar$ and taking imaginary part, we get
\[\hbar\Im(\psi^\dg\nabla\psi)=\hbar\nabla\theta +\hbar\Im(\tilde{\psi}^\dg\nabla\tilde{\psi}),\]
whose divergence vanishes. Consequently, $\psi =e^{\i\theta}\tilde{\psi}$ satisfies the divergence-free condition. Furthermore, it still gives the same director field $n$ since
\[\psi^\dg\sigma_k\psi=e^{-\i\theta}\tilde{\psi}^\dg\sigma_ke^{\i\theta}\tilde{\psi} = \tilde{\psi}^\dg\sigma_k\tilde{\psi}=n_k.\]

In the following, we verify that this constructed spinor $\psi$ indeed satisfies the desired ISF equation \eqref{eq:ISF}.

\subsection{Verification of \eqref{eq:ISF}}
In order to show that $\psi=e^{\i\theta}U(n,n_0)\psi_0$ is a solution to \eqref{eq:ISF}, we need a basic lemma from linear algebra.

\begin{lemma}\label{lem:linalg}
	Suppose $\psi\in\bbc^2$ satisfies $|\psi|^2=1$. Then, for $\phi\in \bbc^2$, $\phi=\i \alpha\psi$ for some $\alpha\in\bbr$ if and only if $\Re(\psi^\dg\sigma \phi)=0$.
\end{lemma}

\begin{proof}
	Consider the linear map $L:\bbc^2\simeq \bbr^4\to\bbr^3$ defined as $L(\phi):=\Re(\psi^\dg\sigma\phi)$. Then, for any real $a\in\bbr^3$, we use \eqref{Pauli_identity_2} and $\psi^\dg\sigma\psi=n\in\bbr^3$ to obtain
	\begin{align*}
		L((a\cdot\sigma)\psi)_k= \Re\left(\psi^\dg \sigma_k(a\cdot\sigma)\psi\right) =\Re\left(\psi^\dg\left(a_kI+\i(e_k\times a)\cdot\sigma\right)\psi\right)=a_k.
	\end{align*}
	Therefore, $L$ is onto, and hence, $\dim\ker L=1$. Observe that for any $\alpha\in\bbr$,
	\[L(\i\alpha\psi)=\Re\left(\psi^\dg\sigma (\i\alpha\psi)\right)=\Re\left(\i\alpha n\right)=0.\]
	Thus, we conclude that
	\[\ker L=\{\phi\in\bbc^2~|~\phi=\i\alpha\psi\quad\mbox{for some}\quad\alpha\in\bbr\}.\]
\end{proof}

Now, we are ready to prove that $\psi$ is indeed a solution to \eqref{eq:ISF}.

\begin{proposition}\label{prop:psi}
	Let $\psi_0\in H^3(\Omega;\bbs^3)$ satisfy $\pa_\nu\psi_0=0$ on $\pa\Omega$ and $\dv\,\Im(\psi_0^\dg\nabla\psi_0)=0$ in $\Omega$, let $n_0=\psi_0^\dg\sigma\psi_0$, and let $n$ be the strong solution to \eqref{eq:ISF_director} given by Theorem \ref{thm:director}.  Let $\psi:=e^{\i\theta}U(n,n_0)\psi_0$, where $\theta$ is the mean-zero solution to the Poisson equation \eqref{eq:poisson}. Then, $\psi$ satisfies
	\[\psi\in L^\infty(0,T_{**};H^3),\quad \pa_t\psi\in L^\infty(0,T_{**};H^1),\quad\pa_\nu\psi=0,\quad \psi(0)=\psi_0,\]
	and it is a strong solution to the ISF equation \eqref{eq:ISF}. Furthermore, $\psi$ can be extended to $T_*$, the existence time of $n$.
\end{proposition}

\begin{proof}
	Below, we show that the constructed $\psi$ satisfies the desired constraints, regularity, and the governing equation \eqref{eq:ISF}.\\
	
	\noindent $\bullet$ (Incompressibility, boundary and initial conditions):	First of all, thanks to the choice of $\theta$, we already observed that
	\[\Im(\psi^\dg\Delta \psi) = \dv\,\Im(\psi^\dg\nabla\psi)=0,\]
	and
	\[\pa_\nu\psi = \pa_\nu (e^{\i\theta}\tilde{\psi})=e^{\i\theta}\pa_\nu\theta\tilde{\psi} +e^{\i\theta}\pa_\nu\tilde{\psi}=0.\]
	Moreover, since $U(n_0,n_0)=I$, we have $\tilde{\psi}(0)=\psi_0$. Furthermore, $\theta(0,x)$ satisfies
	\[\hbar\Delta\theta(0,x)=-\dv (\hbar\Im(\psi^\dg_0\nabla\psi_0))=0,\quad \int_\Omega \theta(0,x)\,dx=0,\]
	which implies $\theta(0,x)\equiv0$ for all $x\in\Omega$. Therefore,
	\[\psi(0,x)=e^{\i\theta(0,x)}\tilde{\psi}(0,x)=\psi_0(x),\quad x\in\Omega.\]
	
	\noindent $\bullet$ (Regularity of $\psi$ and $p$): Since $\tilde{\psi}=U(n,n_0)\psi_0$ is a multiplications of the smooth functions $n_0,n\in L^\infty H^3$ and $\psi_0\in H^3$, it also belongs to $L^\infty H^3$. Furthermore, since 
	\[\pa_t\tilde{\psi}=(\pa_n U\cdot\pa_t n)\psi_0,\] 
	which belongs to $L^\infty H^1$, thanks to the regularity of $\pa_t n$. Furthermore, since the right-hand side \eqref{eq:poisson} belongs to $L^\infty H^1$, we have $\theta\in L^\infty H^3$, and by taking time derivative on \eqref{eq:poisson}, we also get $\pa_t\theta\in L^\infty(0,T_{**};H^1)$. Therefore, we can conclude that $\psi=e^{\i\theta}\tilde{\psi}$ has the desired regularity. Finally, since $p=\psi^\dg(\i\hbar\pa_t\psi+\frac{\hbar^2}{2}\Delta\psi)$, it belongs to $L^\infty H^1$ thanks to the regularity condition $\psi\in L^\infty H^3\subset L^\infty W^{1,\infty}$.\\

	\noindent $\bullet$ (Governing equation): Now, to show that $\psi$ satisfies the ISF equation \eqref{eq:ISF}, we need to show that the residual $R:=\pa_t\psi -\frac{\i\hbar}{2}\Delta\psi$ is purely imaginary multiple of $\psi$, which is $-\frac{\i}{\hbar}p\psi$. Using Lemma \ref{lem:linalg}, it is equivalent to show that $\Re(\psi^\dg\sigma R)=0$. Observe that
	\begin{align*}
		2\Re(\psi^\dg\sigma R)&=2\Re\left(\psi^\dg\sigma \left(\pa_t\psi-\frac{\i\hbar}{2}\Delta\psi\right)\right)=2\Re(\psi^\dg\sigma\pa_t\psi)+\hbar\Im(\psi^\dg\sigma\Delta\psi)\\
		&=\pa_tn +\hbar\pa_j\Im\left(\psi^\dg\sigma\pa_j\psi\right),
	\end{align*}	
	where we used $\Im(\pa_j\psi^\dg\sigma \pa_j\psi)=0$. On the other hand, we already computed in \eqref{Imw} that
	\[\Im\left(\psi^\dg\sigma \pa_j\psi\right)=\frac{v_j}{\hbar}n +\frac{1}{2}\pa_j n\times n,\]
	where $v=\hbar\Im(\psi^\dg\nabla\psi)$ denotes the velocity field generated by $\psi$. Since $\psi$ is a unit spinor satisfying $\dv\,v=0$ and $\pa_\nu\psi=0$, the representation \eqref{represent_v}, whose derivation in Section \ref{sec:2} used only these two properties and not the equation \eqref{eq:ISF}, shows that $v=K[n]$, where $n=\psi^\dg\sigma\psi$ is the solution to \eqref{eq:ISF_director} we started from. Therefore, using $\dv\,v=0$, we obtain
	\[2\Re(\psi^\dg\sigma R)=\pa_tn +\pa_j \left(v_jn +\frac{\hbar}{2}\pa_j n\times n\right)=\pa_tn + (v\cdot \nabla)n-\frac{\hbar}{2} n\times \Delta n=0,\]
	since $n$ solves \eqref{eq:ISF_director}. This proves that $\psi$ satisfies the equation \eqref{eq:ISF}.
	
	To show that $\psi$ can be extended to $T_*$, we first note that $n\in C([0,T_*];H^{3-\delta})$ and $|n(t,x)|^2=1$ for all $t\in[0,T_*]$ and $x\in\Omega$. Therefore, there exists a partition
	\[0=t_0<t_1<\cdots<t_N=T_*\]
	such that $n(t,x)\cdot n(t_i,x)>-1$ for $t\in [t_i,t_{i+1}]$ and all $x\in \Omega$. Hence, for given initial data $\psi_0$, we can inductively construct $\psi$ over all $t\in[0,T_*]$ as
	\[\psi(t,x)=e^{\i\theta(t,x)}U(n(t,x),n(t_i,x))\psi(t_i,x),\quad t\in[t_i,t_{i+1}],\quad x\in\Omega.\]
	After the same verification as above on each subinterval, one can show that this extended $\psi$ satisfies $\psi\in L^\infty(0,T_*;H^3)$, $\pa_t\psi\in L^\infty(0,T_*;H^1)$, $\pa_\nu\psi=0$, and \eqref{eq:ISF}.
\end{proof}

\subsection{Uniqueness of the solution}

As we mentioned in Remark \ref{rem:uniqueness}, the solution $\psi$ to \eqref{eq:ISF} should be normalized so that $\int_\Omega p(t,x)\,dx=0$ in order to guarantee the uniqueness. Therefore, for the solution $(\psi,p)$ constructed in Proposition \ref{prop:psi}, we replace $\psi$ by $e^{\i\vartheta(t)}\psi$ with
\[\vartheta(t):=\frac{1}{\hbar|\Omega|}\int_0^t \int_\Omega p(s,x)\,dx\,ds.\]
By Remark \ref{rem:uniqueness}, $e^{\i\vartheta}\psi$ is again a strong solution in the sense of Definition \ref{def:spinor} with the same initial data, and the pressure $p$ is shifted to $p-\hbar\vartheta'(t)=p-\frac{1}{|\Omega|}\int_\Omega p\,dx$, satisfying the normalization condition. In conclusion, we have constructed the strong solution $(\psi,p)$ to \eqref{eq:ISF} satisfying
\begin{equation}\label{eq:p_regular}
	\int_\Omega p(t,x)\,dx = 0\quad\mbox{for}\quad t\in[0,T_*].
\end{equation}
We now prove that the strong solution $(\psi,p)$ to \eqref{eq:ISF} in the sense of Definition \ref{def:spinor} satisfying \eqref{eq:p_regular} is unique. Suppose $(\psi,p)$ and $(\tilde{\psi},\tilde{p})$ are both strong solutions to \eqref{eq:ISF} on $[0,T_*]$ with the common initial data $\psi_0$, satisfying \eqref{eq:p_regular}, and let $n=\psi^\dg\sigma\psi$ and $\tilde{n}=\tilde{\psi}^\dg\sigma\tilde{\psi}$. Then, following the computation in Section \ref{sec:2}, which consists of pointwise identities between functions having the regularity of Definition \ref{def:spinor}, both $n$ and $\tilde{n}$ are strong solutions to \eqref{eq:ISF_director} in the sense of Definition \ref{def:director} with the common initial data $n_0:=\psi_0^\dg\sigma\psi_0$. Therefore, by the uniqueness of the solution to \eqref{eq:ISF_director}, they should coincide, i.e., $n\equiv \tilde{n}$. Consequently, we also obtain
\[\hbar\Im (\psi^\dg\nabla\psi) = v =K[n]=K[\tilde{n}]=\tilde{v}=\hbar\Im(\tilde{\psi}^\dg\nabla\tilde{\psi}),\]
and
\[\psi\psi^\dg = \frac{1}{2}(I+n\cdot\sigma)=\frac{1}{2}(I+\tilde{n}\cdot\sigma)=\tilde{\psi}\tilde{\psi}^\dg.\]
Now, we define $\zeta:=\psi^\dg\tilde{\psi}\in\bbc$. Then,
\[\tilde{\psi}=\tilde{\psi}(\tilde{\psi}^\dg\tilde{\psi})=(\tilde{\psi}\tilde{\psi}^\dg)\tilde{\psi}=(\psi\psi^\dg)\tilde{\psi}=\zeta\psi.\]
Since $|\psi|=|\tilde{\psi}|=1$, we have $|\zeta|=1$. Moreover, $\zeta(t,\cdot)=\psi^\dg\tilde{\psi}\in H^3(\Omega)\subset C^1(\bar{\Omega})$ and $\Omega$ is simply connected, so that $\zeta=e^{\i\eta}$ for some real-valued $\eta(t,\cdot)\in C^1(\bar{\Omega})$, that is $\tilde{\psi}=e^{\i\eta}\psi$. However, this implies
\[\tilde{v}(t,x)= v(t,x)+\hbar\nabla\eta(t,x).\]
Therefore $\nabla\eta(t,x)=0$, that is $\eta=\eta(t)$ is a function of $t$ only with $\eta(0)=0$. Then, since $\tilde{\psi} = e^{\i\eta}\psi$, which implies 
\[\tilde{\psi}^\dg \pa_t\tilde{\psi} = e^{-\i\eta}\psi^\dg(e^{\i\eta}\pa_t\psi+\i\eta'e^{i\eta}\psi)=\psi^\dg\pa_t\psi+\i\eta',\quad \tilde{\psi}^\dg\Delta\tilde{\psi} = e^{-\i\eta}\psi^\dg \Delta(e^{\i\eta}\psi) = \psi^\dg\Delta\psi,\]
the pressures $p$ and $\tilde{p}$, which are given by $p=\psi^\dg(\i\hbar\pa_t\psi+\frac{\hbar^2}{2}\Delta\psi)$ and $\tilde{p}=\tilde{\psi}^\dg(\i\hbar\pa_t\tilde{\psi}+\frac{\hbar^2}{2}\Delta\tilde{\psi})$ by \eqref{eq:ISF} and $|\psi|=|\tilde{\psi}|=1$, satisfy
\[\tilde{p}=\tilde{\psi}^\dg\left(\i\hbar\pa_t\tilde{\psi}+\frac{\hbar^2}{2}\Delta\tilde{\psi}\right)=\i\hbar(\psi^\dg\pa_t\psi +\i\eta')+\frac{\hbar^2}{2}\psi^\dg\Delta\psi=p-\hbar\eta'.\]
Since we require that the solution is normalized as $\int_\Omega p\,dx=\int_\Omega \tilde{p}\,dx=0$, this implies $\eta'(t)=0$ for $t\in[0,T_*]$. Hence, $\eta$ is constant, and $\eta(0)=0$ gives $\eta\equiv0$. This completes the proof of the uniqueness and thereby, the proof of Theorem \ref{thm:spinor}.

\end{document}